\documentclass[aap,MSNbibl,seceqn,dvips]{arximspdf}
\usepackage{algorithm, algorithmicx, algpseudocode}
\usepackage{graphicx}
\doi{10.1214/13-AAP943} 
\volume{24}
\issue{3}
\pubyear{2014}
\firstpage{1129}
\lastpage{1171}

\makeatletter
\newcommand{\rrvert}{\vert}
\newcommand{\llvert}{\vert}
\newcommand{\eqref}[1]{(\ref{#1})}
\newcommand{\rset}{\mathbb{R}}
\newcommand{\nset}{\mathbb{N}}
\newcommand{\ind}{\mathbf{1}}
\newcommand{\fl}{\longrightarrow}
\newcommand{\e}{\mathbb{E}}
\newcommand{\V}{\mathbb{V}}
\newcommand{\p}{\mathbb{P}}
\newcommand{\D}{\mathbb{D}}
\newcommand{\lp}{\mathrm{L}}
\newcommand{\spp}{\mathrm{S}}
\newcommand{\hpp}{\mathrm{H}}
\newcommand{\m}{\mathcal}
\newcommand{\cf}{\mathcal{F}}
\newcommand{\cc}{\mathcal{C}}

\def\bY{\mathbf{Y}}
\def\bZ{\mathbf{Z}}
\def\bG{\mathbf{G}}
\newtheorem{thmm}{Theorem}[section]
\newtheorem{lemme}[thmm]{Lemma}

\newtheorem{prop}[thmm]{Proposition}

\newtheorem{hypo}[thmm]{Hypothesis}

\newproclaim{df}{Definition}
\newproclaim{ex}[thmm]{Example}
\newproclaim{rem}[thmm]{Remark}
\makeatother

\begin{document}
\begin{frontmatter}

\title{Simulation of BSDEs by Wiener chaos expansion}
\runtitle{Simulation of BSDEs by Wiener chaos expansion}

\begin{aug}
\author[A]{\fnms{Philippe} \snm{Briand}\ead[label=e1]{philippe.briand@univ-savoie.fr}}
\and
\author[A]{\fnms{C\'{e}line} \snm{Labart}\corref{}\ead[label=e2]{celine.labart@univ-savoie.fr}}
\runauthor{P. Briand and C. Labart}
\affiliation{Universit\'{e} de Savoie}
\address[A]{Laboratoire de Math\'{e}matiques\\
CNRS UMR 5127\\
Universit\'{e} de Savoie\\
Campus Scientifique\\
73376 Le Bourget du Lac\\
France\\
\printead{e1}\\
\hphantom{E-mail:\ }\printead*{e2}} 
\end{aug}

\received{\smonth{4} \syear{2012}}
\revised{\smonth{3} \syear{2013}}

%
\begin{abstract} We present an algorithm to solve BSDEs based on Wiener chaos
expansion and Picard's iterations. We get a forward scheme where the
conditional expectations are easily computed thanks to chaos decomposition
formulas. We use the Malliavin derivative to compute $Z$. Concerning the
error, we derive explicit bounds with respect to the number of chaos
and the
discretization time step. We also present numerical experiments. We obtain
very encouraging results in terms of speed and accuracy.
\end{abstract}

%
\begin{keyword}[class=AMS]
\kwd{60H10}
\kwd{60H35}
\kwd{65C05}
\kwd{65G99}
\kwd{60H07}
\end{keyword}
\begin{keyword}
\kwd{Backward stochastic differential equation}
\kwd{Wiener chaos expansion}
\end{keyword}

\end{frontmatter}

\section{Introduction} 
\label{secintroduction} In this paper, we are interested in the numerical
approximation of solutions $(Y,Z)$ to backward stochastic differential
equations (BSDEs for short in the sequel). BSDEs were introduced by
Bismut in \cite{Bis73} in the linear case, whereas the nonlinear case
was considered later by Pardoux and Peng in \cite{PP90}. A~BSDE is
an equation of the following form:
%
\begin{equation}
\label{eqmain} Y_t=\xi+\int_t^T
f(s,Y_s,Z_s) \,ds - \int_t^T
Z_s\cdot dB_s, \qquad 0\leq t\leq T,
\end{equation}
where $B$ is a $d$-dimensional standard Brownian motion, the terminal
condition $\xi$ is a real-valued $\m F_T$-measurable random variable where
$\{\m F_t\}_{0\leq t\leq T}$ stands for the augmented filtration of the
Brownian motion $B$ and the generator $f$ is a map from $[0,T] \times
\rset\times\rset^d$ into $\rset$. A solution to this equation is a
pair of
processes $\{(Y_t,Z_t)\}_{0\leq t\leq T}$ which is required to be
adapted to
the filtration $\{\m F_t\}_{0\leq t\leq T}$. We will assume the conditions
of Pardoux and Peng to ensure existence and uniqueness of solutions.

Our main objective in this study is the numerical approximation of the
solution $(Y,Z)$ to BSDE~\eqref{eqmain} (even though there exists a large
literature on this subject). The first two contributions to this topic
are due
to Chevance~\cite{Che97}, who considered generators independent of $Z$
and Bally~\cite{Bal97}, who used a random time mesh. Ma and Yong
\cite{MY99} proposed numerical schemes based on the link between Markovian
\mbox{BSDEs} and semilinear partial differential equations (PDEs). Another approach,
based on Donsker's theorem and close to \cite{Che97}, was proposed by
Coquet, Mackevicius and M\'emin \cite{CMM99} in the case of a
generator $f$ independent of $Z$; the general case was treated by Briand,
Delyon and M\'emin in \cite{BDM01}, who later extended it to a more
general framework \cite{BDM02}, including the case of
a ``stepwise constant Brownian motion.'' This extension led to the formulas
\[
Y_t = \e (Y_{t+h} | \m F_t ) + h
f(t,Y_t,Z_t),\qquad Z_t = h^{-1} \e
\bigl(Y_{t+h} (B_{t+h}-B_t ) | \m F_t
\bigr)
\]
known as the dynamic programming algorithm. Even though the convergence was
proved in the case of path-dependent terminal condition $\xi$, the
rate of
convergence was left as an open question in \cite{BDM02}. This problem was
solved by Zhang \cite{Zha04} and Bouchard and Touzi \cite{BT04} in
the case of Markovian BSDE, namely in the case of a terminal condition
$\xi=g(X_T)$ where $X$ is the solution to a stochastic differential
equation; in \cite{Zha04}, the author considers the path-dependent case
as well. Their result was generalized by Gobet and Labart \cite{GL07}
and also by Gobet and Makhlouf~\cite{GM10}.

From a numerical point of view, the main difficulty in solving BSDEs
is to
efficiently compute conditional expectations. Several approaches have been
proposed using various tools: the Malliavin calculus \cite{BT04}, regression
methods \cite{GLW05,GLW06} and quantization technics \cite{BP03}.


Finally, let us mention that there exist some works dealing with the
discretization of solutions to BSDEs in a more general framework:
forward--backward SDEs~\cite{DM06} and quadratic BSDEs \cite{Ric11}.

Let us now describe briefly the main points of our approach in the case
of a real-valued Brownian motion. Already used in several quoted papers
(see also \cite{BD07,GL10,BS11}), our starting point is the use of
Picard's iterations, $(Y^0,Z^0)=(0,0)$ and for $q\in\nset$,
\[
Y^{q+1}_t = \xi+ \int_t^T
f \bigl(s,Y^q_s,Z^q_s \bigr)
\,ds - \int_t^T Z^{q+1}_s
\cdot dB_s,\qquad 0\leq t\leq T.
\]
It is well known that the sequence $(Y^q,Z^q)$ converges exponentially
fast toward the solution $(Y,Z)$ to BSDE~\eqref{eqmain}.
We write this Picard scheme in a forward way,
\begin{eqnarray*}
Y^{q+1}_t & = &\e \biggl(\xi+ \int_0^T
f \bigl(s,Y^q_s,Z^q_s \bigr)
\,ds\Big | \m F_t \biggr) - \int_0^t f
\bigl(s,Y^q_s,Z^q_s \bigr) \,ds,
\\
Z^{q+1}_t & =& D_t Y^{q+1}_t
= D_t \e \biggl(\xi+ \int_0^T f
\bigl(s,Y^q_s,Z^q_s \bigr) \,ds
\Big| \m F_t \biggr),
\end{eqnarray*}
where $D_t X$ stands for the Malliavin derivative of the random
variable $X$.

In order to compute the previous conditional expectation, we use a
Wiener chaos expansion of the random variable
\[
F^{q}=\xi+ \int_0^T f
\bigl(s,Y^q_s,Z^q_s \bigr) \,ds.
\]
More precisely, we use the following orthogonal decomposition of the
random variable $F^q$:
\[
F^q = \e \bigl[F^q \bigr] + \sum
_{k\geq1} \sum_{|n| = k}
d_k^n \prod_{i\geq1}
K_{n_i} \biggl(\int_0^T
g_i(s)\,dB_s \biggr),
\]
where $K_l$ denotes the Hermite polynomial of degree $l$, $(g_i)_{i\geq
1}$ is an
orthonormal basis of $\lp^2(0,T)$ and, if $n=(n_i)_{i\geq1}$ is a sequence
of integers, $|n| = \sum_{i\geq1} n_i$. $(d^n_k)_{k \ge1, |n|=k}$ is the
sequence of coefficients ensuing from the decomposition of $F^q$. Of course,
from a practical point of view, we only keep a finite number of terms
in this expansion:
\begin{itemize}
\item we work with a finite number of chaos, $p$;
\item we choose a finite number of functions $g_1, \ldots, g_N$.
\end{itemize}
This leads to the following approximation with $n=(n_1,\ldots,n_N)$:
\[
F^q \simeq\e \bigl[F^q \bigr] + \sum
_{1\leq k\leq p} \sum_{|n| =
k}
d_k^n \prod_{1\leq i\leq N}
K_{n_i} \biggl(\int_0^T
g_i(s)\,dB_s \biggr).
\]
One of the key points in using such a decomposition is that, for
choices of
simple functions $g_1$, \ldots, $g_N$, there exist explicit formulas
for both
%
\begin{equation}
\label{eq30} \e \bigl( F^q | \m F_t \bigr) \quad\mbox{and}\quad
Z^{q+1}_t = D_t \e \bigl( F^q | \m
F_t \bigr);
\end{equation}
this plays a crucial role in our algorithm. Using these formulas and
starting from $M$ trajectories of the underlying Brownian motion, we
are able to construct $M$ trajectories of the solution $(Y,Z)$ to the BSDE.

In the following, the functions $g_i$ are chosen as step functions:
\[
g_i=\ind_{]\overline{t}_{i-1}, \overline{t}_i]}(t)/ \sqrt{h},\qquad i=1,\ldots,N,
 \mbox{ where }
\overline{t}_i:=ih, h=\frac{T}{N}
\]
and the previous formulas are really simple; see
\eqref{chaosdec}--\eqref{coefchaosdec} and
Proposition~\ref{prop2}. Eventually, the
main advantage of this method is that only one decomposition has
to be computed per Picard iteration: the decomposition of $F^q$.
Therein lies
the main difference between our approach and the approach based on
regression technique developed by Bender and Denk in \cite{BD07}. In
their paper, for a given Picard iteration $q$ and for each time $t_i$
of the mesh
grid, two projections have to be computed, one for $Y^q_{t_i}$ and one
for $Z^q_{t_i}$. The second
difference comes from the way of computing $Z^q$. In our method, once
the decomposition of $F^q$ is computed, $Z^q$ is given explicitly as
the Malliavin derivative of $Y^q$. Let us also point out that our
algorithm can handle fully path dependent terminal
conditions.


The rest of the paper is organized as follows. Section~\ref{secwienerchaosexpansion} contains the notation and the preliminary
results, Section~\ref{secdescalgo} describes precisely the algorithm,
Section~\ref{secconvres} is devoted to the study of the convergence of
the algorithm and finally Section~\ref{secnumex} contains some
numerical experiments. Some technical proofs are postponed to the
\hyperref[app]{Appendix}.\vadjust{\goodbreak}

%
\section{Preliminaries} 
\label{secwienerchaosexpansion}
%

\subsection{Definitions and notation}
Given a probability space $(\Omega,\cf,\p)$ and an $\rset^d$-valued Brownian
motion $B$, we consider:
\begin{itemize}
\item$\{(\cf_t);t\in[0,T]\}$, the filtration generated by the
Brownian motion
$B$ and augmented.
\item$\lp^p(\m F_T):=\lp^p(\Omega,\cf_T,\p)$, $p\in\nset^*$, the
space of all $\cf_T$-measurable random variables
(r.v. in the following) $X\dvtx \Omega\longmapsto\rset^d$ satisfying
$\|X\|^p_p: = \e(|X|^p)< \infty$.
\item$\e_t(X)$ denotes $\e(X |\cf_t)$ for any $X$ in $\lp^1(\cf_T)$.

\item$\spp^p_T(\rset^d)$, $p \in\nset, p\ge2$, the space of all
c\`
{a}dl\`{a}g predictable processes $\phi\dvtx \Omega\times[0,T]
\longmapsto
\rset^d$ such that $\|\phi\|^p_{\spp^p}=\e
(\sup_{t \in[0,T]} |\phi_t|^p)< \infty$.
\item$\hpp^p_T(\rset^d)$, $p \in\nset, p\ge2$, the space of all
predictable processes
$\phi\dvtx \Omega\times[0,T] \longmapsto\rset^d$ such that $\|\phi
\|^p_{\hpp^p_T}=\e
\int_0^T |\phi_t|^p \,dt < \infty$.
\item$\lp^2(0,T)$, the space of all square integrable functions on $[0,T]$.
\item$C^{k,l}$, the set of continuously differentiable functions $\phi\dvtx (t,x)
\in[0,T] \times\rset^d$ with continuous derivatives w.r.t. $t$
(resp., w.r.t. $x$) up to order $k$ (resp., up to order $l$).
\item$C^{k,l}_b$, the set of continuously differentiable functions
$\phi\dvtx (t,x) \in[0,T] \times\rset^d$ with continuous and uniformly bounded
derivatives w.r.t. $t$ (resp., w.r.t. $x$) up to order~$k$ (resp., up to order
$l$). The function $\phi$ is also bounded.
\item$\|\partial^j_{\mathrm{sp}} f\|^2_{\infty}$, the norm of the
derivatives of
$f([0,T] \times\rset^d, \rset)$ w.r.t. all the space
variables $x$ which sum equals $j$: $\|\partial^j_{\mathrm{sp}}
f\|^2_{\infty}:=\sum_{|k|=j}
\|\partial^{k_1}_{x_1} \cdots\partial^{k_d}_{x_d} f\|^2_{\infty
}$, where
$|k|=k_1+\cdots+k_d$.
%
\item$C^{\infty}_p$, the set of smooth functions $f\dvtx \rset^n
\longmapsto\rset$
with partial derivatives of polynomial growth.
\item$\|(\cdot,\cdot)\|^p_{\lp^p}$, $p \in\nset, p\ge2$, the norm
on the space $\spp^p_T(\rset)\times
\hpp^p_T(\rset^d)$ defined by
%
\begin{equation}
\label{normeL2YZ}\bigl \|(Y,Z)\bigr\|^p_{\lp^p}:= \e\Bigl( \sup
_{t \in[0,T]} |Y_t|^p\Bigr)+\int
_0^T \e\bigl(|Z_t|^p
\bigr) \,dt.
\end{equation}
\end{itemize}
We also recall some useful definitions related to Malliavin
calculus. We use the notation of \cite{nualart06}.
\begin{itemize}
\item$\mathcal{S}$ denotes the class of random variables of the form
$F=f(W(h_1),\ldots,\break W(h_n))$, where $f \in C^{\infty}_p(\rset
^{n\times
d},\rset)$, for all $j \le n$
$h_j=(h_j^1,\ldots,h_j^d) \in\lp^2([0,T];\break \rset^d)$ and for all
$i\le
d$ $W^i(h^i_j)=\int_0^T
h_j^i(t)\,dW^i_t$.
\item$\mathbb{D}^{r,2}$ denotes the closure of $\mathcal{S}$ w.r.t.
the following norm on $\mathcal{S}$
\[
\|F\|^2_{\mathbb{D}^{r,2}}:=\e|F|^2 + \sum
_{q=1}^r \sum_{|\alpha
|_1=q}\e
\biggl(\int_0^T \cdots\int_0^T
\bigl\llvert D^{\alpha}_{(t_1,\ldots,t_q)}F\bigr\rrvert ^2
\,dt_1\cdots dt_q \biggr),
\]
where $\alpha$ is a multi-index $(\alpha_1,\ldots,\alpha_q) \in
\{1,\ldots,d\}^q$ $|\alpha|_1:=\sum_{i=1}^q \alpha_i=q$, and
$D^{\alpha
}$ represents the
multi-index Malliavin derivative operator. We recall $\mathbb
{D}^{\infty,2}=\bigcap_{r=1}^{\infty}
\mathbb{D}^{r,2}$.
\end{itemize}

\begin{rem}
When $d=1$, $\|F\|^2_{\mathbb{D}^{r,2}}:=\e|F|^2 +\sum_{q=1}^r
\e
(\int_0^T \cdots\break \int_0^T
\llvert D^{(q)}_{(t_1,\ldots,t_q)}F\rrvert ^2 \,dt_1\cdots dt_q
)=\e
|F|^2 +\sum_{q=1}^r \|D^{(q)}F\|^2_{\lp^2(\Omega\times[0,T]^q)}$.
\end{rem}
Let $m \in\nset^*$ and $j \in\nset,j\ge2$. We also introduce the
following notation:
\begin{itemize}
\item$\mathcal{D}^{m,j}$ denotes the
space of all $\cf_T$-measurable r.v. such that
\[
\|F\|^j_{m,j}:=\sum_{1\le l \le m}
\sum_{|\alpha|_1=l} \sup_{t_1\le
\cdots\le t_l} \e
\bigl[\bigl|D^{\alpha}_{t_1,\ldots,t_l} F\bigr|^j\bigr]< \infty,
\]
where $\sup_{t_1\le\cdots\le t_l}$ means $\sup_{(t_1,\ldots,t_l) \dvtx t_1\le\cdots\le t_l}$.
\item$\mathcal{S}^{m,j}$ denotes the space of all couple of processes $(Y,Z)$
belonging to $\spp^j_T(\rset)\times
\hpp^j_T(\rset^d)$ and such that
\[
\bigl\|(Y,Z)\bigr\|^j_{m,j}:=\sum_{1\le l \le m}
\sum_{|\alpha|_1=l} \sup_{t_1\le\cdots\le t_l}\bigl\|
\bigl(D^{\alpha}_{t_1, \ldots, t_l} Y,D^{\alpha}_{t_1, \ldots, t_l} Z\bigr)
\bigr\|^j_{\lp^j}<\infty.
\]
We recall
\begin{eqnarray*}
&&\bigl\|(Y,Z)\bigr\|^j_{m,j}= \sum_{1\le l \le m}
\sum_{|\alpha|_1=l} \sup_{t_1\le\cdots\le t_l} \biggl\{ \e
\Bigl[\sup_{t_l\le r \le T} \bigl|D^{\alpha}_{t_1, \ldots, t_l}
Y_r\bigr|^j\Bigr] \\
&&\hspace*{152pt}{}+\int_{t_l}^T
\e\bigl[\bigl|D^{\alpha}_{t_1, \ldots, t_l} Z_r\bigr|^j\bigr]
\,dr \biggr\}.
\end{eqnarray*}
We also denote $\mathcal{S}^{m,\infty}:=\bigcap_{j \ge2} \mathcal{S}^{m,j}$.
\end{itemize}

\subsection{Wiener chaos expansion}
\subsubsection{Notation and useful results}
We refer to \cite{nualart06} for more
details on this section. Let us briefly recall the Wiener chaos
expansion in the simple case of a real-valued Brownian\vadjust{\goodbreak} motion. It is
well known that every
random variable $F\in\lp^2(\m F_T)$ has an expansion of the following form:
%
\begin{eqnarray}
\label{eq1} F &=& \e[F]+ \int_0^T
u_1(s_1)\,dB_{s_1}
\nonumber
\\
&&{}+ \int_0^T \int_0^{s_2}
u_2(s_2,s_1)\,dB_{s_1}\,dB_{s_2}
+ \cdots\\
&&{}+ \int_0^T \int_0^{s_{n}}
\cdots\int_{0}^{s_2} u_n(s_n,
\ldots,s_1) \,dB_{s_1}\cdots dB_{s_n} + \cdots,\nonumber
\end{eqnarray}
where the functions $(u_n, n \ge1)$ are deterministic functions. There
is an ambiguity
for the definition of these functions $u_n$. We adopt in this paper the
following point of view: the function $u_n$ is defined on the simplex
\[
\mathcal{S}_n(T):= \bigl\{(s_1,\ldots,s_n)
\in[0,T]^n \dvtx 0<s_1<\cdots <s_n<T \bigr\}.
\]
We define the iterated integral for a
deterministic function $f \in\lp^2(\mathcal{S}_n(T))$ as
\[
J_n(f):= \int_0^T \int
_0^{s_{n}} \cdots\int_0^{s_2}
f(s_n,\ldots,s_1) \,dB_{s_1} \cdots
dB_{s_n}.
\]
Due to the It\^{o} isometry,
$\|J_n(f)\|^2=\|f\|^2_{\lp^2(\mathcal{S}_n(T))}$ and
$\e[J_n(f) J_m(g)]=\break \delta_{nm} \langle f,g\rangle _{\lp^2(\mathcal{S}_n(T))}$. Then
$\|F\|^2=\sum_{n \ge0}
\|u_n\|^2_{\lp^2(\mathcal{S}_n(T))}$.

\begin{df*}
Let $F$ be a random variable in $\lp^2(\m F_T)$ whose chaos expansion
is given by \eqref{eq1}. We introduce:
\begin{itemize}
\item$P_n(F):= J_n({u_n})$ the Wiener chaos of order $n$ of $F$.
\item$\cc_p(F):=\sum_{n \le p} P_n(F)$ the chaos decomposition of $F$
up to
order $p$, that is,
%
\begin{eqnarray}
\label{cp} \cc_p(F)&=&\e[F]+ \int_0^T
u_1(s_1)\,dB_{s_1}+ \int_0^T
\int_0^{s_2} u_2(s_2,s_1)\,dB_{s_1}\,dB_{s_2}
\nonumber
\\[-8pt]
\\[-8pt]
\nonumber
&&{}+ \cdots+ \int_0^T \int_0^{s_{p}}
\cdots\int_{0}^{s_2} u_p(s_p,
\ldots,s_1) \,dB_{s_1}\cdots dB_{s_p}.
\end{eqnarray}
\end{itemize}
\end{df*}

We state two lemmas useful for the sequel.


\begin{lemme}[(Nualart)]\label{lem1}
$F \in\D^{m,2}$ if and only if $\|D^m F\|^2_{\lp^2(\Omega\times
[0,T]^m)}=\sum_{n \ge0} (n+m-1)\times\cdots
\times n \times\e[|P_n(F)|^2] < \infty$. In this case, we have
\[
\sum_{n \ge0} (n+m-1)\times\cdots \times n \times\e
\bigl[\bigl|P_n(F)\bigr|^2\bigr] \leq\|F\|^2_{\mathbb{D}^{m,2}}.
\]
\end{lemme}

From Lemma~\ref{lem1}, we deduce the following:
%
\begin{lemme}\label{lem2}
Let $F \in\D^{m,2}$. We have
\[
\e\bigl[\bigl|F-\cc_p(F)\bigr|^2\bigr] \le\frac{\|D^m F\|^2_{\lp^2(\Omega\times
[0,T]^m)}}{(p+m)\cdots(p+1)}.
\]
\end{lemme}

\begin{pf}
\begin{eqnarray*}
\e\bigl[\bigl|F-\cc_p(F)\bigr|^2\bigr] &=&\sum
_{k \ge p+1} \e\bigl[P_k(F)^2\bigr]\\
&=&\sum
_{k \ge p+1} (k+m-1)\cdots k\times\frac{1}{(k+m-1)\cdots k}\times\e
\bigl[\bigl|P_k(F)\bigr|^2\bigr]
\\
&\le&\frac{1}{(p+m)\cdots(p+1)}\sum_{k \ge p+1} (k+m-1)\cdots k \e
\bigl[\bigl|P_k(F)\bigr|^2\bigr].
\end{eqnarray*}
\upqed\end{pf}

The following lemma gives some useful properties of the chaos decomposition.
%
\begin{lemme}\label{lem4}
\begin{itemize}
\item Let $F$ be a r.v. in $\lp^2(\m F_T)$. $\forall p \ge1$, we have $
\e(|\cc_p(F)|^2) \le\e(|F|^2)$. If $F$ belongs to $\lp^j(\m F_T)$,
$\forall j>2$,
$\e(|\cc_p(F)|^j) \le(1+p(j-1)^{{p}/{2}})^j \e(|F|^j)$.
\item Let $H$ be in $\hpp^2_T(\rset)$. We have $\cc_p  (\int_0^T
H_s \,ds )=\int_0^T \cc_p (H_s) \,ds$.
\item For all $F \in\mathbb{D}^{1,2}$ and for all $t \le r$, $D_t \e
_r[\cc_p(F)]=\e_r[\cc_{p-1}(D_t F)]$.
\end{itemize}
\end{lemme}
The first result ensues from the fact that for $j>2$ $\|P_n(F)\|_j \le
(j-1)^{{n}/{2}}\|F\|_j$; see \cite{nualart06}, page~63.

\subsubsection{Wiener chaos expansion and Hermite polynomials} Another
approach to Wiener chaos expansion uses Hermite polynomials. This approach
can be easily generalized when considering $d$-dimensional Brownian
motions, and so this is the one we consider in the following. We
present it
for $d=1$. Let $\{g_i\}_{i\geq1}$ be an orthonormal basis
of $\lp^2(0,T)$. The Wiener chaos of order $n$, $ P_n(F)$, is the
$\lp^2$-closure of the vector field spanned by
\[
\biggl\{ \prod_{i\geq1} \sqrt{n_i!}
K_{n_i} \biggl(\int_0^T
g_i(s) \,dB_s \biggr) \dvtx \bigl|(n_i)_{i\geq1}\bigr|:=
\sum n_i =n \biggr\},
\]
where $K_n$ is the Hermite polynomial of order $n$ defined by the expansion
\[
e^{xt-t^2/2} = \sum_{n\geq0} K_n(x)
t^n
\]
%
with the convention $K_{-1}\equiv0$. With this normalization, we have
$K_n'(x) = K_{n-1}(x)$ for any integer $n$.
It is well known that $(K_n)_{n\geq0}$ is a sequence of orthogonal
polynomials in
$\lp^2(\rset,\mu)$,\vadjust{\goodbreak} where $\mu$ denotes the reduced centered Gaussian
measure. Moreover, we have
\[
\int_{\rset} K_n^2(x) \mu(dx) =
\frac{1}{n!}.
\]

Every square integrable random variable $F$, measurable with respect to
$\m
F_T$, admits the following orthogonal decomposition:
%
\begin{equation}
\label{eqtruedec} F = d_0 + \sum_{k\geq1}
\sum_{|n| = k} \,d_k^n \prod
_{i\geq
1} K_{n_i} \biggl(\int
_0^T g_i(s)\,dB_s
\biggr),
\end{equation}
where $n=(n_i)_{i\geq1}$ is a sequence of positive integers, and where
$|n|$ stands for $\sum_{i\geq1} n_i$. Taking into account the
normalization of the Hermite polynomials we use, we get
\[
d_0 = \e [F ],\qquad d_k^n = n! \e \biggl[F
\times\prod_{i\geq1} K_{n_i} \biggl(\int
_0^T g_i(s)\,dB_s
\biggr) \biggr],
\]
where $n! = \prod_{i\geq1} n_i !$. Before describing the chaos decomposition
formulas we use in the algorithm, we give a lemma useful in the sequel.

\begin{lemme}\label{enmart}
Let $g\in\lp^2(0,T)$, and let $U_t = \int_0^t g^2(s) \,ds$. For $n\in
\nset$, let us define
\[
M^n_t = U_t^{n/2} K_n
\bigl( B(g)_t / \sqrt{U_t} \bigr),\qquad B(g)_t =
\int_0^t g(s)\,dB_s.
\]

Then $\{M^n_t\}_{0\leq t\leq T}$ is a martingale and
\[
dM^n_t = g(t) M^{n-1}_t
\,dB_t.
\]
\end{lemme}

\subsection{Chaos decomposition formulas}
These formulas are based on the decomposition \eqref{eqtruedec}. To get
tractable formulas, we consider a finite number of chaos and a finite number
of functions $(g_1,\ldots,g_N)$. The $(g_i)_{1 \le i \le N}$ functions are
chosen such that we can quickly compute $\e(F| \m F_t)$ and $D_t \e(F|
\m
F_t)$ [as required in \eqref{eq30}]. We develop in this section the
case $d=1$, and
we refer to Section~\ref{subbmRd} when $d>1$.

The first step consists in considering a finite number of chaos. In
order to
approximate the random variable $F$, we consider its projection $\cc_p(F)$
onto the first $p$ chaos, namely
%
\begin{equation}
\label{eqcp} \cc_p(F) = d_0 + \sum
_{1\leq k\leq p} \sum_{|n| = k}
d_k^n \prod_{i\geq1}
K_{n_i} \biggl(\int_0^T
g_i(s)\,dB_s \biggr).
\end{equation}

Of course, we still have an infinite number of terms in the previous
sum and
the second step consists in working with only the first $N$ functions
$g_1,\ldots, g_N$ of an orthonormal basis of $\lp^2(0,T)$.

Let us consider a regular mesh grid of $N$ time steps
$\mathcal{T}=\{\overline{t}_i=i\frac{T}{N}, i=0,\ldots,N \}$ and the
$N$ step
functions
%
\begin{equation}
\label{eqgi} g_i=\ind_{]\overline{t}_{i-1}, \overline{t}_i]}(t)/ \sqrt{h},\qquad i=1,\ldots,N,
\mbox{ where } h:=\frac{T}{N}.
\end{equation}
We complete these $N$ functions $g_1,\ldots, g_N$ into an
orthonormal basis
of $\lp^2(0,T)$, $(g_i)_{i\geq1}$. For instance, one can consider the Haar
basis on each interval $(\overline{t}_{i-1},\overline{t}_i)$,
$i=1,\ldots, N$. We implicitly assume that $N\geq p$. This leads to
the following approximation:
%
\begin{equation}
\label{eqcpN} \cc_p^N(F) = d_0 + \sum
_{1\leq k\leq p} \sum_{|n| = k}
d_k^n \prod_{1\leq i\leq N}
K_{n_i} \biggl(\int_0^T
g_i(s)\,dB_s \biggr),
\end{equation}
where $n=(n_1,\ldots,n_N)$ and $|n|=n_1+\cdots+n_N$. Due to the
simplicity of the functions $g_i$, $i=1,\ldots,N$, we can compute explicitly
\[
\int_0^T g_i(s)\,dB_s
= G_i\qquad \mbox{where } G_i = \frac{B_{\overline
{t}_i}-B_{\overline{t}_{i-1}}}{\sqrt{h}}.
\]
Roughly speaking this means that $P_k$, the $k$th chaos, is generated by
\[
\bigl\{ K_{n_1}(G_1)\cdots K_{n_N}(G_N)
\dvtx n_1+\cdots+n_N = k \bigr\}.
\]
Thus the approximation we will use for the random variable $F$ is
%
\begin{eqnarray}
\label{chaosdec} \cc_p^N(F) &=& d_0 + \sum
_{k=1}^p \sum
_{|n|=k} \,d_k^n K_{n_1}(G_1)
\cdots K_{n_N}(G_N)
\nonumber
\\[-8pt]
\\[-8pt]
\nonumber
&=& d_0 + \sum
_{k=1}^p \sum_{|n|=k}
d_k^n \prod_{1\leq
i\leq
N}
K_{n_i}(G_i),
\end{eqnarray}
where the coefficients $d_0$ and $d_k^n$ are given by
%
\begin{equation}
\label{coefchaosdec} d_0=\e[F], \qquad d_k^n = n! \e
\bigl[F K_{n_1}(G_1)\cdots K_{n_N}(G_N)
\bigr].
\end{equation}

The following lemma, similar to Lemma~\ref{lem4}, gives
some useful properties of the operator $\cc^N_p$.
%
\begin{lemme}\label{lem9}
Let $F$ be a r.v. in $\lp^2(\m F_T)$ and $H$ be in $\hpp^2_T(\rset
)$. Then:
\begin{itemize}
\item$\forall(p,N) \in(\nset^{\star})^2, \e(|\cc^N_p(F)|^2) \le
\e(|\cc_p(F)|^2)\le\e(|F|^2)$.
\item$\cc^N_p  (\int_0^T H_s \,ds )=\int_0^T \cc^N_p
(H_s) \,ds$.
\item For all $t \le r$, $D_t \e_r[\cc^N_p(F)]=\e_r[\cc^N_{p-1}(D_t F)]$.
\end{itemize}
\end{lemme}

From \eqref{chaosdec}, we deduce the expressions of $\e_t(\cc_p^N
F)$ and
$D_t\e_t ( \cc_p^N(F) )$, useful for the approximation of
$(Y,Z)$ by
the chaos decomposition; see \eqref{eq30}.

\begin{prop}\label{prop2} Let $F$ be a real random variable in $\lp
^2(\m F_T)$, and let $r$ be an integer in $\{1,\ldots, N\}$.
For all $\overline{t}_{r-1}<t\leq\overline{t}_r$, we have
\begin{eqnarray*}
\e_t \bigl(\cc_p^N F \bigr)& =&
d_0 + \sum_{k=1}^p \sum
_{|n(r)|=k} d_k^n \prod
_{i<r} K_{n_i}(G_i)\\
&&\hspace*{19pt}{}\times \biggl(
\frac{t-\overline
{t}_{r-1}}{h} \biggr)^{{n_r}/{2}} K_{n_r} \biggl(
\frac
{B_t-B_{\overline
{t}_{r-1}}}{\sqrt{t-\overline{t}_{r-1}}} \biggr),
\\
D_t \e_t \bigl(\cc_p^N(F)
\bigr) 
&=& h^{-1/2} \sum_{k=1}^p
\mathop{\sum_{|n(r)|=k
}}_{ n_r>0}
d_k^n \prod_{i<r}
K_{n_i}(G_i)\\
&&{}\times \biggl(\frac{t-\overline{t}_{r-1}}{h}
\biggr)^{{(n_r-1)}/{2}} K_{n_r-1}
\biggl(\frac{B_t-B_{\overline{t}_{r-1}}}{\sqrt{t-\overline
{t}_{r-1}}} \biggr),
\end{eqnarray*}
where, if $r\leq N$ and $n=(n_1,\ldots,n_N)$, $n(r)$ stands for
$(n_1,\ldots,n_r)$.
\end{prop}

The proof of Proposition~\ref{prop2} is postponed to Section~\ref{subbmR}.

\begin{rem}\label{rem1}
For $t=\overline{t}_r$ and $r\ge1$, Proposition~\ref{prop2} leads to
\begin{eqnarray*}
\e_{\overline{t}_r} \bigl(\cc_p^N F \bigr) & =&
d_0 + \sum_{k=1}^p \sum
_{|n(r)|=k} d_k^n \prod
_{i\leq r} K_{n_i} (G_i ),
\\
D_{\overline{t}_r}\e_{\overline{t}_r} \bigl(\cc_p^N
F \bigr) & =& h^{-1/2} \sum_{k=1}^p
\mathop{\sum_{|n(r)|=k }}_{ n_r>0}
d_k^n \prod_{i<r}
K_{n_i} (G_i )\times K_{n_r-1} (G_r
).
\end{eqnarray*}
When $r=0$, we get $\e_{\overline{t}_0} (\cc_p^N F ) =
d_0$,  and
we define $D_{\overline{t}_0}\e_{\overline{t}_0} (\cc_p^N
F )=\frac{1}{\sqrt{h}}d^{e_1}_1$ [which is the limit of
$D_{t}\e_{t} (\cc_p^N F )$ when $t$ tends to $0$].
\end{rem}

Let us end this subsection by some examples.

\begin{ex}[(Case $p=2$)]
From \eqref{chaosdec}--\eqref{coefchaosdec}, we have
\[
\cc^N_2(F)=d_0+\sum
_{j=1}^N d_1^{e_j}
K_1(G_j)+\sum_{j=1}^N
\sum_{i=1}^{j-1}d_2^{e_{ij}}K_1(G_i)K_1(G_j)+
\sum_{j=1}^N d_2^{2e_j}
K_2(G_j),
\]
where $e_j$ denotes the unit vector whose $j$th component is one, and
$e_{ij}=e_i+e_j$. For $j=1,\ldots,N$ and $i=1,\ldots,j-1$, it holds
\begin{eqnarray*}
d_1^{e_j}&=&\e\bigl(FK_1(G_j)\bigr),\qquad
d_2^{e_{ij}}=\e\bigl(FK_1(G_i)K_1(G_j)
\bigr), \\
 d_2^{2e_j}&=&2\e\bigl(FK_2(G_j)
\bigr).
\end{eqnarray*}
Remark~\ref{rem1} leads to
\begin{eqnarray*}
\e_{\overline{t}_r} \bigl(\cc^N_2 F \bigr)
&=&d_0+\sum_{j=1}^r
d_1^{e_j} K_1(G_j)+\sum
_{j=1}^r\sum_{i=1}^{j-1}d_2^{e_{ij}}K_1(G_i)K_1(G_j)\\
&&{}+
\sum_{j=1}^r d_2^{2e_j}
K_2(G_j),
\\
 D_{\overline{t}_r} \e_{\overline{t}_r} \bigl( \cc_2^N F
\bigr) &=&h^{-1/2} \Biggl(d_1^{e_r}+d^{2e_r}_2K_1(G_r)+
\sum_{i=1}^{r-1}d_2^{e_{ir}}
K_1(G_i) \Biggr).
\end{eqnarray*}
\end{ex}

%

\section{Description of the algorithm}
\label{secdescalgo}
The algorithm is based on four types of approximations: Picard's iterations,
a Wiener
chaos expansion up to a finite order, the truncation of an $\lp^2(0,T)$
basis in order to apply formulas of Proposition~\ref{prop2}, and a
Monte Carlo
method to approximate the coefficients $d_0$ and $d^n_k$ defined in
\eqref{coefchaosdec}. We present the
first three steps of the approximation procedure in Section~\ref{sectalgotheorique}. The Monte Carlo method and the practical
implementation are presented in Section~\ref{sectimplementation}.



\subsection{Approximation procedure}\label{sectalgotheorique}
\subsubsection{Picard's iterations}\label{sectpicard}
The first step consists in approximating $(Y,Z)$---the
solution to \eqref{eqmain}---by Picard's sequence
$(Y^q,Z^q)_q$, built as follows:
$(Y^0=0,Z^0=0)$ and for all $q \ge1$
%
\begin{equation}
\label{eqpicard} Y^{q+1}_t = \xi+\int_t^T
f \bigl(s,Y^q_s,Z^q_s \bigr)
\,ds - \int_t^T Z^{q+1}_s
\cdot dB_s, \qquad 0\leq t\leq T.
\end{equation}

From \eqref{eqpicard}, under the assumptions that $\xi\in\D^{1,2}$
and $f \in C^{0,1,1}_b$, we
express $(Y^{q+1},Z^{q+1})$ as a function of the processes $(Y^q,Z^q)$,
%
\begin{equation}
\label{eqYqZq} Y^{q+1}_t= \e_t \biggl( \xi+ \int
_t^T f \bigl(s,Y^q_s,Z^q_s
\bigr) \,ds \biggr),\qquad Z^{q+1}_t=D_t
Y^{q+1}_t,
\end{equation}
which can also be written
%
\begin{eqnarray}
\label{eqpY} Y^{q+1}_t &=& \e_t \biggl( \xi+
\int_0^T f \bigl(s,Y^q_s,Z^q_s
\bigr) \,ds \biggr) - \int_0^t f
\bigl(s,Y^q_s,Z^q_s \bigr) \,ds,
\nonumber
\\[-8pt]
\\[-8pt]
\nonumber
Z^{q+1}_t&=&D_t Y^{q+1}_t.
\end{eqnarray}

As we recalled in the \hyperref[secintroduction]{Introduction}, the computation of the
conditional expectation is the cornerstone in the numerical resolution of
BSDEs. Chaos decomposition formulas enable us to circumvent this
problem.

\subsubsection{Wiener Chaos expansion} Computing the chaos
decomposition of the
r.v. $F=\xi+ \int_t^T f (s,Y^q_s,Z^q_s ) \,ds$ [appearing in
\eqref{eqYqZq}] in order to compute $Y^{q+1}_t$ is not judicious. $F$ depends
on $t$, and then the computation of $Y^{q+1}$ on the grid
$\mathcal{T}=\{\overline{t}_i=i\frac{T}{N}, i=0,\ldots,N \}$ would
require $N+1$ calls to
the chaos decomposition function. To build an efficient algorithm, we
need to
call the chaos decomposition function as infrequently as possible,
since each call
is computationally demanding and brings an approximation error due to the
truncation and to the Monte Carlo approximation (see next sections).
Then we
look for a r.v. $F^q$ independent of $t$ such that $Y^{q+1}_t$ and $Z^{q+1}_t$
can be expressed as functions of $\e_t(F^q)$, $D_t \e_t (F^q)$ and of
$Y^q$ and
$Z^q$. Equation \eqref{eqpY} gives a more tractable expression of
$Y^{q+1}$. Let $F^q$ be defined by $F^q:= \xi+ \int_0^T f(s,Y^q_s,Z^q_s)
\,ds$. Then
%
\begin{equation}
\label{eqY} Y^{q+1}_t= \e_t
\bigl(F^q\bigr)- \int_0^t f
\bigl(s,Y^q_s,Z^q_s \bigr) \,ds,\qquad
Z^{q+1}_t=D_t \e_t
\bigl(F^q\bigr).
\end{equation}

The second type of approximation consists of computing the chaos decomposition
of $F^q$ up to order $p$. Since $F^q$ does not depend on $t$, the chaos
decomposition function $\cc_p$ is called only once per Picard's iteration.

Let $(Y^{q,p},Z^{q,p})$ denote the approximation of $(Y^q,Z^q)$ built
at step
$q$ using a chaos decomposition with order $p$:
$(Y^{0,p},Z^{0,p})=(0,0)$ and
%
\begin{eqnarray}
\label{eqYpZp} Y^{q+1,p}_t&=&\e_t \bigl[
\cc_p \bigl(F^{q,p} \bigr) \bigr]-\int_0^t
f \bigl(s,Y^{q,p}_s,Z^{q,p}_s
\bigr)\,ds,
\nonumber
\\[-8pt]
\\[-8pt]
\nonumber
Z^{q+1,p}_t&=&D_t \e_t \bigl[
\cc_p \bigl(F^{q,p} \bigr) \bigr],
\end{eqnarray}
where $F^{q,p}=\xi+\int_0^T
f (s,Y^{q,p}_s,Z^{q,p}_s )\,ds$.
In the sequel, we also use the following equality:
%
\begin{equation}
\label{eqZp} Z^{q+1,p}_t= \e_t
\bigl[D_t \cc_p\bigl(F^{q,p}\bigr)\bigr].
\end{equation}


\subsubsection{Truncation of the basis} The third type of
approximation comes
from the truncation of the orthonormal $\lp^2(0,T)$ basis used in the
definition of $\cc_p$~\eqref{eqcp}. Instead of
considering a basis of $\lp^2(0,T)$, we only keep the first $N$ functions
$(g_1,\ldots,g_N)$ defined by \eqref{eqgi} to build the chaos decomposition
function $\cc_p^N$ \eqref{eqcpN}. Proposition~\ref{prop2} gives us explicit formulas for $\e_t (\cc_p^N
F)$ and $D_t \e_t (\cc_p^N
F)$. From~\eqref{eqYpZp}, we build $((Y^{q,p,N},Z^{q,p,N})_q$ in the following
way: $((Y^{0,p,N},Z^{0,p,N})=(0,0)$ and
%
\begin{eqnarray}
\label{eqYZqpN} Y^{q+1,p,N}_t&=& \e_t\bigl(
\cc^N_p \bigl(F^{q,p,N}\bigr)\bigr)- \int
_0^t f \bigl(s,Y^{q,p,N}_s,Z^{q,p,N}_s
\bigr) \,ds,
\nonumber
\\[-8pt]
\\[-8pt]
\nonumber
 Z^{q+1,p,N}_t&=&D_t\bigl(\e _t
\bigl(\cc^N_p \bigl(F^{q,p,N}\bigr)\bigr)\bigr),
\end{eqnarray}
where $F^{q,p,N}:= \xi+ \int_0^T f(s,Y^{q,p,N}_s,Z^{q,p,N}_s)
\,ds$.\vadjust{\goodbreak}

Equation \eqref{eqYZqpN} is tractable as soon as we know closed
formulas for the\allowbreak coefficients $d^n_k$ of the chaos decomposition of
$\e_t(\cc^N_p (F^{q,p,N}))$ and\break $D_t(\e_t(\cc^N_p (F^{q,p,N})))$; see
Proposition~\ref{prop2}. When it is not the case, we need to use a
Monte Carlo
method to approximate these coefficients. The next section is devoted
to this
method and to the
practical implementation. In particular, we give the pseudo-code of the
algorithm.

\subsection{Implementation}\label{sectimplementation}
In this section, we first explain how to practically compute the chaos
decomposition $\cc_p^N(F)$ of a r.v. $F$. Then we give the pseudo-code
of the algorithm.

\subsubsection{Monte Carlo simulations of the chaos
decomposition}\label{secMC} Let $F$ denote a r.v. of
$\lp^2(\m F_T)$. Practically, when we are not able to compute exactly
$d_0$ and/or
the coefficients $d^n_k$ of the chaos decomposition
\eqref{chaosdec}--\eqref{coefchaosdec} of $F$, we use Monte Carlo
simulations to approximate them. Let $(F^m)_{1\le m \le M}$ be a $M$
i.i.d. sample of $F$ and $(G_1^m,\ldots,G_N^m)_{1\le m \le M}$ be a $M$
i.i.d. sample of $(G_1,\ldots,G_N)$. We recall that $d_0$ and the
coefficients $(d^n_k)_{1
\le k \le p,|n|=k}$ are given by $ d_0=\e[F]$ and $d_k^n = n! \e [F
K_{n_1}(G_1)\cdots K_{n_N}(G_N) ]$; see \eqref{coefchaosdec}. Then
they are solutions of
%
\begin{equation}
\label{LS} \mathop{\arg\min}_{\mathbf{c}=(c_0,(c^n_k)_{1 \le k \le p,|n|=k})}
 \e \bigl[\bigl|F-
\psi(c,G)\bigr|^2\bigr],
\end{equation}
where $\psi\dvtx (\mathbf{c},G) \longmapsto c_0 + \sum_{k=1}^p \sum_{|n|=k}
c_k^n \prod_{1\leq i\leq N} K_{n_i}(G_i)$. We propose two methods to
approximate $\mathbf{d}:=(d_0,(d^n_k)_{1 \le k \le p,|n|=k})$:

\begin{itemize}
\item the first one consists in approximating the expectations of
\eqref{coefchaosdec} by empirical
means $ \mathbf{\widehat{d_M}}:=(\hat{d_0},\hat{d^n_k}_{1 \le k \le
p,|n|=k})$ where
%
\begin{equation}
\label{dchapeau} \widehat{d_0}:=\frac{1}{M}\sum
_{m=1}^M F^m,\qquad \widehat{d^n_k}:=
\frac
{n!}{M}\sum_{m=1}^M
F^m K_{n_1}\bigl(G_1^m\bigr)
\cdots K_{n_N}\bigl(G_N^m\bigr);
\end{equation}
\item the second one is based on a sample average approximation
\[
\mathbf{\overline{d_M}}:=\bigl(\overline{d_0},
\overline{d^n_k}_{1 \le k
\le p,|n|=k}\bigr)= \mathop{\arg
\min}_{c_0,(c^n_k)_{1 \le k \le
p,|n|=k}} \frac{1}{M}\sum_{m=1}^M
\bigl|F^m-\psi\bigl(\mathbf{c},G^m\bigr)\bigr|^2.
\]
\end{itemize}

\begin{rem}\label{rem4}
In terms of computation time, the first method is much faster than the second
one.
\begin{itemize}
\item The first method requires
$O(M\times p)$ computations per coefficient. Since we are looking for
$O(N^p)$ coefficients, its computational cost is
$O(M\times p \times N^p)$.
\item The second method requires $O(M\times p \times N^p)$ computations to
evaluate $\frac{1}{M}\sum_{m=1}^M |F^m-\psi(c,G^m)|^2$ (in fact, it requires
the same number of computations as the first method, since the function
$\psi$ contains as many additions as coefficients, and each addition
contains as many products as the associated coefficient). We still have
to compute the argmin, the computational cost of which depends on the
method we
use.
\end{itemize}
%

From a theoretical point of view, the second method gives better convergence
results than the first one. For the first method, we only know that
$\mathbf{\widehat{d_M}}$ converges
to~$\mathbf{d}$ a.s. Concerning the second method, we know that~$\mathbf
{\overline{d_M}}$ converges to~$\mathbf{d}$ a.s., and under regularity
assumptions on $\psi$, the uniform strong law of large numbers gives
the a.s.
convergence of $\frac{1}{M}\sum_{m=1}^M
|F^m-\psi(\mathbf{\overline{d_M}},G^m)|^2$ to
$\e[|F-\psi(\mathbf{d},G)|^2]$.
\end{rem}
In the following, $\cc_p^{N,M} (F)$ denotes the
approximation of the chaos decomposition of order $p$ of $F$ when using
the first method to approximate the coefficients~$d^n_k$:
%
\begin{equation}
\label{chaosdecMC} \cc_p^{N,M} (F)=\widehat{d_0} +
\sum_{k=1}^p \sum
_{|n|=k} \widehat {d_k^n} \prod
_{1\leq i\leq N} K_{n_i}(G_i).
\end{equation}
$\e_t(\cc_p^{N,M} (F))$ and
$D_t(\e_t(\cc_p^{N,M} (F)))$ denote the conditional expectations
obtained in
Proposition~\ref{prop2} when $(d_0,d^n_k)_{1 \le k \le p,|n|=k}$ are
replaced by $(\widehat{d_0},\break \widehat{d^n_k})_{1 \le k \le p,|n|=k}$,
\begin{eqnarray*}
\e_t \bigl(\cc_p^{N,M} F \bigr)&:=&
\widehat{d_0} + \sum_{k=1}^p
\sum_{|n(r)|=k} \widehat{d_k^n}
\prod_{i<r} K_{n_i}(G_i)\\
&&\hspace*{19pt}{}\times
\biggl(\frac{t-\overline{t}_{r-1}}{h} \biggr)^{{n_r}/{2}} K_{n_r} \biggl(
\frac{B_t-B_{\overline{t}_{r-1}}}{\sqrt{t-\overline
{t}_{r-1}}} \biggr),
\\
D_t \e_t \bigl(\cc_p^{N,M}(F)
\bigr) 
&:=& h^{-1/2} \sum_{k=1}^p
\mathop{\sum_{|n(r)|=k
}}_{ n_r>0}
\widehat{d_k^n} \prod_{i<r}
K_{n_i}(G_i)\\
&&{}\times \biggl(\frac{t-\overline{t}_{r-1}}{h}
\biggr)^{{(n_r-1)}/{2}} K_{n_r-1}
\biggl(\frac{B_t-B_{\overline{t}_{r-1}}}{\sqrt{t-\overline
{t}_{r-1}}} \biggr).
\end{eqnarray*}

\begin{rem}\label{rem2} When $M$ samples of $\cc_p^{N,M} (F)$ are
needed, we
can either use the same samples as the ones used to compute $\widehat{d_0}$
and $\widehat{d_k^n}$: $(\widehat{\cc_p^N} (F))^m=\widehat{d_0} +
\sum_{k=1}^p \sum_{|n|=k} \widehat{d_k^n} \prod_{1\leq i\leq N}
K_{n_i}(G^m_i)$, or use new ones. In the first case, we only require $M$
samples of $F$ and $(G_1,\ldots,G_N)$. The coefficients $\widehat{d_k^n}$
and $\widehat{d_0}$ are not independent of $\prod_{1\leq i\leq N}
K_{n_i}(G^m_i)$. The notation $\e_t(\cc_p^{N,M} (F))$ introduced
above cannot
be linked to $\e (\cc_p^{N,M} F|\cf_t )$. In the second
case, the
coefficients $\widehat{d_k^n}$ and $\widehat{d_0}$ are independent of
$\prod_{1\leq i\leq N} K_{n_i}(G^m_i)$, and we have $\e_t (\cc_p^{N,M}
F )=\e (\cc_p^{N,M} F|\cf_t )$. This second approach
requires $2M$
samples of $F$ and $(G_1,\ldots,G_N)$, and its variance increases with
$N$. Practically, we use the first technique.
\end{rem}

We introduce the processes $(Y^{q+1,p,N,M},Z^{q+1,p,N,M})$, which is
useful in the
following. It corresponds to the
approximation of $(Y^{q+1,p,N},Z^{q+1,p,N})$ when we use $\cc_p^{N,M}$ instead
of $\cc_p^N$, that is, when we use a Monte Carlo procedure to compute
the coefficients
$d^n_k$.
%
\begin{eqnarray}
\label{eqYZqpNM} Y^{q+1,p,N,M}_t&=& \e_t\bigl(
\cc^{N,M}_p \bigl(F^{q,p,N,M}\bigr)\bigr)- \int
_0^t f \bigl({\theta}^{q,p,N,M}_s
\bigr) \,ds, Z^{q+1,p,N,M}_t
\nonumber
\\[-8pt]
\\[-8pt]
\nonumber
&=&D_t\bigl(\e _t
\bigl(\cc^{N,M}_p \bigl(F^{q,p,N,M}\bigr)\bigr)\bigr),
\end{eqnarray}
where $F^{q,p,N,M}:= \xi+ \int_0^T
f({\theta}^{q,p,N,M}_s) \,ds$ and ${\theta}^{q,p,N,M}_s=
(s,{Y}^{q,p,N,M}_s,\allowbreak {Z}^{q,p,N,M}_s )$.

\subsubsection{Pseudo-code of the algorithm}\label{sectionalgo} In
this section, we describe in
details the algorithm. We aim at computing $M$
trajectories of an approximation of $(Y,Z)$ on the grid
$\mathcal{T}=\{\overline{t}_i=i\frac{T}{N}, i=0,\ldots,N \}$.
Starting from
$(Y^{0,p,N,M},Z^{0,p,N,M})=(0,0)$, \eqref{eqYZqpNM} enables to get
$(Y^{q,p,N,M},Z^{q,p,N,M})$ for each of Picard's iterations $q$ on
$\mathcal{T}$. Practically, we discretize the integral $ \int_0^t
f ({\theta}^{q,p,N,M}_s  ) \,ds$ which leads to approximated
values of
$(Y^{q,p,N,M},Z^{q,p,N,M})$ computed on a grid.

Let us introduce
$(\overline{Y}^{q+1,p,N,M}_{\overline{t}_i},\overline
{Z}^{q+1,p,N,M}_{\overline{t}_i})_{1\le i \le
N}$, defined by $(\overline{Y}^{0,p,N,M},\allowbreak \overline
{Z}^{0,p,N,M})=(0,0)$ and
for all $q \ge0$
%
\begin{eqnarray}
\label{eqYZqpNMbarre}\quad\qquad \overline{Y}^{q+1,p,N,M}_{\overline{t}_i}&=& \e_{\overline{t}_i}
\bigl(\cc ^{N,M}_p \bigl(\overline{F}^{q,p,N,M}\bigr)
\bigr)-h \sum_{j=1}^{i} f \bigl(
\overline{t}_j,\overline{Y}^{q,p,N,M}_{\overline
{t}_j},\overline
{Z}^{q,p,N,M}_{\overline{t}_j} \bigr),
\nonumber
\\[-8pt]
\\[-8pt]
\nonumber
\overline{Z}^{q+1,p,N,M}_{\overline{t}_i}&=&D_{\overline{t}_i}\bigl(\e
_{\overline{t}_i}\bigl(\cc^{N,M}_p\bigl(\overline{F}^{q,p,N,M}
\bigr)\bigr)\bigr),
\end{eqnarray}
where $\overline{F}^{q,p,N,M}:= \xi+ h\sum_{i=1}^{N}
f(\overline{t}_i,\overline{Y}^{q,p,N,M}_{\overline{t}_i},\overline
{Z}^{q,p,N,M}_{\overline{t}_i})$.
Here is the notation we use in the algorithm:

\begin{itemize}
\item$d$: dimension of the Brownian motion;
\item$q$: index of Picard's iteration;
\item$K_{it}$: number of Picard's iterations;
\item$M$: number of Monte Carlo samples;
\item$N$: number of time steps used for the discretization of $Y$ and $Z$;
\item$p$: order of the chaos decomposition;
\item$\bY^q \in\mathcal{M}_{N+1,M}(\rset)$ represents $M$ paths of
$\overline{Y}^{q,p,N,M}$
computed on the grid~$\mathcal{T}$;
\item for all $l \in\{1,\ldots,d\}$, $(\bZ^q)_l \in
\mathcal{M}_{N+1,M}(\rset)$ represents $M$ paths of $(\overline
{Z}^{q,p,N,M})_l$ computed on
the grid $\mathcal{T}$.
\end{itemize}

Since $\xi\in\lp^2(\m F_T)$, $\xi$ can be written as a measurable function
of the Brownian path. Then one gets one sample of $\xi$ from one
sample of
$(G_1,\ldots,G_N)$ (where $G_i$ represents $\frac{B_{\overline
{t}_i}-B_{\overline{t}_{i-1}}}{\sqrt{h}}$).

For the sake of clarity, we detail the algorithm for $d=1$.
\begin{algorithm}[t]
\caption{Iterative algorithm}\label{algodetail}
\begin{algorithmic}[1]
\State Pick at random $N\times M$ values of standard Gaussian
r.v. stored in $\bG$.
\State Using $\bG$, compute $(\xi^m)_{0 \le m \le M-1}$.
\State$\bY^0 \equiv0$, $\bZ^0 \equiv0$.
\For{$q=0\dvtx K_{it}-1$} \label{Kit}
\For{$m=0\dvtx M-1$}\label{loopF}
\State Compute $(F^q)^m=\xi^m+h\sum_{i=1}^{N}f(\overline{t}_i,(\bY
^q)_{i,m},(\bZ^q)_{i,m})$
\EndFor
\State Compute the vector $\mathbf{d}=(\widehat{d_0},(\widehat
{d^n_k})_{1\le
k \le p, |n|=k})$ of the chaos decomposition of $F^q$\label{coef}
\State
$\widehat{d_0}:=\frac{1}{M}\sum_{m=0}^{M-1}
(F^q)^m, \widehat{d^n_k}=\frac{n!}{M}\sum_{m=0}^{M-1} (F^q)^m
K_{n_1}(G_1^m)\cdots
K_{n_N}(G_N^m)$
\For{$j=1\dvtx N$}\label{lineN}
\For{$m=0\dvtx M-1$}\label{lineM}
\State Compute
$(\e_{\overline{t}_j}(\cc_p^{N,M} F^q))^m$,
$(D_{\overline{t}_j}(\e_{\overline{t}_j}(\cc_p^{N,M}
F^q)))^m$\label{loopEt}
\State$(\bY^{q+1})_{j,m}=(\e_{\overline{t}_j}(\cc_p^{N,M}
F^q))^m-h\sum_{i=1}^{j}f(\overline{t}_i,(\bY^q)_{i,m},(\bZ
^q)_{i,m})$\label{loopYZ}
\State
$(\bZ^{q+1})_{j,m}=(D_{\overline{t}_j}(\e_{\overline{t}_j}(\cc
_p^{N,M} F^q)))^m$
\EndFor
\EndFor
\EndFor
\State Return $(\bY^{K_{it}})_{0,:}=\hat{d}_0$ and $(\bZ
^{K_{it}})_{0,:}=\frac{1}{\sqrt{h}}\hat{d}^{e_1}_1$
\end{algorithmic}
\end{algorithm}

Let us now deal with the complexity of the algorithm:

For each $q$:
\begin{itemize}
\item the computation of the vector $F^q$ (loop line \ref{loopF}) requires
$O(M\times N)$ computations;
\item the computation of the vector $\mathbf{d}$ (line \ref{coef})
requires $O(M\times p \times
(N\times d)^p)$ computations [in dimension $d$ we have $O( (N\times
d)^p)$ coefficients, and the computation of each coefficient requires
$O(M\times p)$
computations (see Remark~\ref{rem4})];
\item for each $N$ and $M$ (lines \ref{lineN}--\ref{lineM}):
\begin{itemize}[$-$]
\item[$-$] the computation of $(\e_{\overline{t}_j}(\cc_p^{N,M} F^q))^m$
and of
$(D^l_{\overline{t}_j}(\e_{\overline{t}_j}(\cc_p^{N,M}
F^q)))^m_{1\le
l \le d}$ (line~\ref{loopEt}) requires
$O(d \times p\times(N\times d)^p)$ computations
\item[$-$] the computation of $(\bY^{q+1})_{j,m}$ (loop line \ref{loopYZ})
requires $O(N)$ computations and the computation of
$((\bZ^{q+1})^l_{j,m})_{1\le l \le d}$ requires $O(d)$ computations.
\end{itemize}
\end{itemize}
The complexity of the algorithm is then $O(K_{it}\times M
\times p\times(N \times d)^{p+1})$.

\section{Convergence results}
\label{secconvres}

We aim at bounding the error between $(Y,Z)$---the solution of
\eqref{eqmain}---and $(Y^{q,p,N,M},Z^{q,p,N,M})$ defined by
\eqref{eqYZqpNM}. Before stating the main result of the paper, we introduce
some hypotheses.

In the following, $(t_1,\ldots,t_n)$ and $(s_1,\ldots,s_n)$ denote
two vectors
such that
\[
0\le t_1\le\cdots\le t_n \le T, 0\le s_1\le
\cdots\le s_n \le T \mbox{ and } \forall i, s_i \le
t_i.
\]

\begin{hypo}[(Hypothesis $\mathcal{H}_m$)]\label{hypo3}
Let $m \in\nset^*$. We say that $F$ satisfies Hypothesis $\mathcal
{H}_m$ if
$F$ satisfies the two following hypotheses:
\begin{itemize}
\item$\mathcal{H}_m^1$: $\forall j \ge2$ $F \in\mathcal
{D}^{m,j}$, that is,
$\|F\|^j_{m,j} < \infty$;
\item$\mathcal{H}_m^2$: $\forall j \ge2$, $\forall i \in
\{1,\ldots,m\}$, $\forall l_0 \le i-1$, $\forall l_1 \le m-i$,
$\forall l
\in\{1,\ldots,d\}$ and for all
multi-indices $\alpha_0$ and $\alpha_1$ such that $|\alpha_0|=l_0$ and
$|\alpha_1|=l_1+1$, there exist two positive constants $\beta_F$ and
$k^F_l$ such that
\begin{eqnarray*}
&&\sup_{t_1\le\cdots\le t_{l_0}}\sup_{s_{i+1}\le\cdots\le s_{i+l_1}} \e\bigl[\bigl|
D^{\alpha_0}_{t_1,\ldots,t_{l_0}}\bigl(D^{\alpha
_1}_{t_i,s_{i+1},\ldots
,s_{i+l_1}}
F-D^{\alpha_1}_{s_i,\ldots,s_{i+l_1}} F \bigr)\bigr|^j\bigr] \\
&&\qquad\le
k^{F}_l(j) (t_i-s_i)^{j \beta_{F}},
\end{eqnarray*}
where $l=l_0+l_1+1$. In the following, we denote $K^F_m(j)=\sup_{l\le
m} k^F_l(j)$.
\end{itemize}
\end{hypo}

\begin{rem}
If $F$ satisfies $\mathcal{H}^2_m$, for all multi-index $\alpha$ such
that $|\alpha|=l$, we have
%
\begin{equation}\quad
\label{eq31} \bigl|\e\bigl(D^{\alpha}_{t_1,\ldots,t_l}F\bigr)-\e
\bigl(D^{\alpha}_{s_1,\ldots
,s_l}F\bigr)\bigr|\le K^F_l
\bigl((t_1-s_1)^{\beta_{F}}+\cdots+(t_l-s_l)^{\beta_{F}}
\bigr),
\end{equation}
where $K^F_l$ is a constant.
\end{rem}

\begin{hypo}[(Hypothesis $\mathcal{H}^3_{p,N}$)]\label{hypo4}
Let $(p,N) \in\nset^2$. We say that an r.v. $F$ satisfies $\mathcal
{H}^3_{p,N}$ if
\[
V_{p,N}(F):=\V(F)+\sum_{k=1}^p
\sum_{|n|=k}n! \V \Biggl(F\prod
_{i=1}^N K_{n_i}(G_i)
\Biggr) < \infty.
\]
\end{hypo}

\begin{rem}\label{rem9}
If $F$ is bounded by $K$, we get $ V_{p,N}(F)\le K^2 \sum_{k=0}^p {N
\choose k}$. Then every bounded r.v. satisfies $\mathcal{H}^3_{p,N}$.
\end{rem}
This remark ensues from $\e (\prod_{i=1}^N
K_{n_i}^2(G_i) )=\frac{1}{n!}$.

\begin{rem}\label{rem6} Let $X$ be the $\rset^n$-valued process
solution of
\[
X_t=x+\int_0^t
b(s,X_s)\,ds +\int_0^t
\sigma(s,X_s)\,dB_s,
\]
where $B$ is a $d$-dimensional Brownian motion and $b\dvtx [0,T] \times
\rset^n
\rightarrow\rset^n$ and $\sigma\dvtx  [0,T] \times\rset^n \rightarrow
\rset^{n\times d}$ are two $C^{0,m}$ functions uniformly Lipschitz
w.r.t. $x$ and H\"{o}lder continuous of parameter $\frac{1}{2}$ w.r.t.
$t$, with linear growth in $x$ and with bounded derivatives.
Then, every random variable $\xi$ of type $g(X_T)$ or $g(\int_0^T X_s \,ds)$
with $g\dvtx \rset^n \rightarrow\rset$ in $C^{\infty}_{p}$ satisfies
$\mathcal{H}_m$ and $\mathcal{H}^3_{p,N}$, for all $p$ and $N$.
\end{rem}
We refer to Section~\ref{proofrem6} for the proof of Remark~\ref{rem6}.

\begin{thmm}\label{thmm1} Let $k$ be an integer s.t. $k \le p$. Assume that
$\xi$ satisfies $\mathcal{H}_{p+q}$ and $\mathcal{H}^3_{p,N}$ and
$f\in
C^{0,p+q-1,p+q-1}_b$. We have
\begin{eqnarray*}
&&\bigl\|\bigl(Y-{Y}^{q,p,N,M},Z-{Z}^{q,p,N,M}\bigr)\bigr\|^2_{\lp
^2}
\\
&&\qquad\le \frac{A_0}{2^q} + \frac{A_1(q,k)}{(p+1)^k}+A_2(q,p) \biggl(
\frac{T}{N} \biggr)^{2\beta_{\xi} \wedge1}+ \frac{A_3(q,p,N)}{M},
\end{eqnarray*}
where $A_0$ is given in Section~\ref{secterreurpicard}, $A_1$
is given in Proposition~\ref{prop1}, $A_2$ is given in Proposition~\ref{prop4}, and $A_3$ is given in Proposition~\ref{prop5}.

If $f\in C^{0,\infty,\infty}_b$ and $\xi$ satisfies $\mathcal
{H}_{\infty}$
and $\mathcal{H}^3_{\infty,\infty}$, we get
\[
\lim_{q \rightarrow\infty}\lim_{p \rightarrow
\infty}\lim
_{N \rightarrow
\infty}\lim_{M \rightarrow
\infty}\bigl\|\bigl(Y-Y^{q,p,N,M},Z-Z^{q,p,N,M}
\bigr)\bigr\|^2_{\lp^2}=0.
\]
\end{thmm}

\begin{rem} If $f$ is a path-dependent generator, Theorem \ref
{thmm1} still
holds true under the following hypotheses: $\forall l \le p$, $\forall j
\ge2$, for all
multi-index $\alpha$ in $\{1,\ldots,d+1\}^l$ ($d$ is the dimension of the
Brownian motion) s.t. $a(i)=d+1$ means that the Malliavin derivative
w.r.t. $t_i$
concerns the path-dependent component, and we assume
\begin{eqnarray*}
&&\int_0^T \bigl\|D^{\alpha}_{t_1,\ldots,t_l}f
\bigl(s,Y^q_s,Z^q_s\bigr)
\bigr\|^2_{\lp
^2(\Omega\times
[0,T]^m)}\,ds < \infty,
\\
&& \int_0^T \e\bigl[\bigl|D^{\alpha}_{t_1,\ldots,t_l}f
\bigl(s,Y^{q}_s,Z^{q}_s
\bigr)\bigr|^j\bigr]\,ds < \infty,\\
&& \int_0^T
\e\bigl[\bigl|D^{\alpha}_{t_1,\ldots,t_l}f\bigl(s,Y^{q,p}_s,Z^{q,p}_s
\bigr)\bigr|^j\bigr]\,ds < \infty\quad\mbox{and }
\\
&&\bigl |\e\bigl(D^{\alpha}_{t_1,\ldots,t_l}I_{q,p}\bigr)-\e
\bigl(D^{\alpha
}_{s_1,\ldots
,s_l}I_{q,p}\bigr)\bigr|\\
&&\qquad\le
K^{I_{q,p}}_l\bigl((t_1-s_1)^{\beta_{I_{q,p}}}+
\cdots+(t_l-s_l)^{\beta
_{I_{q,p}}}\bigr),
\end{eqnarray*}
where $I_{q,p}=\int_0^Tf(\theta_r^{q,p})\,dr$, and $K^{I_{q,p}}_l$ and
$\beta_{I_{q,p}}$ are two positive constants.\vadjust{\goodbreak}
\end{rem}

\begin{rem} Given the complexity $C_0$ of the algorithm (and a given
value of
$d$), we can choose the parameters
$p,q,N$ and $M$ such that they minimize the error $ \frac{A_0}{2^q} +
\frac{A_1(q,p)}{(p+1)^p}+A_2(q,p)  (\frac{T}{N}
)^{a}+\frac{A_3(q,p,N)}{M}$, where $a:=2\beta_{\xi} \wedge1$.
This boilds down to solving the
following constrained minimization problem:
\[
\min_{q,p,N,M\ \mathrm{s.t.}\ qpMN^{p+1}=C_0} \biggl( \frac{1}{2^q} + \frac{C^q}{(p+1)^p}+
\frac{C^q}{N^a}+\frac{C^q N^p}{M} \biggr).
\]
The Karush--Kuhn--Tucker theorem gives $M\sim\frac
{2p}{a}(p+1)^{p+{p^2}/{a}}$, $N\sim(p+1)^{{p}/{a}}$, $q\sim\frac{1}{\ln(2C)}p \ln
(p+1)$ and $p$ such that
$(p+1)^{2p(1+{p}/{a})} p^3 \ln(p+1) \sim a\log(2C)C_0$.
\end{rem}

\begin{pf*}{Proof of Theorem~\ref{thmm1}}
We split the error into $4$ terms:
\begin{longlist}[(1)]
\item[(1)] Picard's
iterations. $\mathcal{E}^q=\|(Y-Y^q,Z-Z^q)\|^2_{\lp^2}$, where
$(Y^q,Z^q)$ is defined by \eqref{eqpicard};
\item[(2)] the truncation of the chaos decomposition.
$\mathcal{E}^{q,p}=\|(Y^q-Y^{q,p},Z^q-Z^{q,p})\|^2_{\lp^2}$, where
$(Y^{q,p},Z^{q,p})$ is defined by \eqref{eqYpZp};
\item[(3)] the
truncation of the $\lp^2(0,T)$ basis.
$\mathcal{E}^{q,p,N}=\|(Y^{q,p}-Y^{q,p,N},Z^{q,p}-Z^{q,p,N})\|^2_{\lp^2}$,
where $(Y^{q,p,N},Z^{q,p,N})$ is defined by \eqref{eqYZqpN};
\item[(4)] the Monte Carlo approximation to compute the
expectations.\break
$\mathcal{E}^{q,p,N,M}=\|
(Y^{q,p,N}-Y^{q,p,N,M},Z^{q,p,N}-Z^{q,p,N,M})\|^2_{\lp^2}$,
where $(Y^{q,p,N,M},\break Z^{q,p,N,M})$ is defined by \eqref{eqYZqpNM}.
\end{longlist}
We have
\[
\bigl\|\bigl(Y-Y^{q,p,N,M},Z-Z^{q,p,N,M}\bigr)\bigr\|^2_{\lp^2}
\le4\bigl(\mathcal {E}^q+\mathcal{E}^{q,p}+
\mathcal{E}^{q,p,N}+\mathcal{E}^{q,p,N,M}\bigr).
\]
It remains to combine \eqref{erreurYZq}, Propositions~\ref{prop1},
\ref{prop4} and~\ref{prop5} to get
the first result.
\end{pf*}

\subsection{Picard's iterations}\label{secterreurpicard}
The first type of error has already been studied in \cite{pardoux92} and
\cite{karoui97}, and we only recall the main result.

\begin{hypo}\label{hypo5}
We assume:
\begin{itemize}
\item the generator $f\dvtx \rset^+\times\rset\times\rset^d\fl\rset$
is Lipschitz
continuous: there exists a constant $L_f$ such that for all $t \in
\rset^+$, $y_1,y_2 \in\rset$ and $z_1,z_2 \in\rset^d$
\[
\bigl|f(t,y_1,z_1)-f(t,y_2,z_2)\bigr| \leq
L_f \bigl(|y_1-y_2|+|z_1-z_2|
\bigr);
\]
\item$\e[|\xi|^2+\int_0^T |f(s,0,0)|^2 \,ds] < \infty$.
\end{itemize}
\end{hypo}
From \cite{karoui97}, Corollary~2.1, we know that under
Hypothesis \ref{hypo5}, the sequence $(Y^q,Z^q)_q$ defined by
\eqref{eqpicard} converges\vadjust{\goodbreak} to $(Y,Z)$ $d \p
\times dt$ a.s. and in $\spp^2_T(\rset)\times\hpp^2_{T}(\rset^d)$.
Moreover, we have
%
\begin{equation}
\label{erreurYZq} \mathcal{E}^q:=\bigl\|\bigl(Y-Y^q,Z-Z^q
\bigr)\bigr\|^2_{\lp^2} \le \frac{A_0}{2^q},
\end{equation}
where $A_0$ depends on $T$, $\|\xi\|^2$
and on $\|f(\cdot,0,0)\|^2_{\lp^2_{(0,T)}}$.

\subsection{Error due to the truncation of the chaos
decomposition}\label{secterreurchaos}
We assume that the integrals are computed exactly, as well as
expectations. The error is only due to the truncation of the chaos
decomposition $\cc_p$ introduced in \eqref{cp}.


For the sequel, we also need the following lemma. We postpone its proof to
the Appendix \ref{prooflem5}.

\begin{lemme}\label{lem5}
Assume that $\xi$ satisfies $\mathcal{H}_{m+q}^1$ and $f \in
\mathcal{C}^{0,m+q-1,m+q-1}_b$. Then $\forall
q' \le q$, $\forall p \in\nset$, $(Y^{q'},Z^{q'})$ and
$(Y^{q',p},Z^{q',p})$ belong
to $\mathcal{S}^{m,\infty}$. Moreover
\begin{eqnarray*}
&&\bigl\|\bigl(Y^{q},Z^{q}\bigr)\bigr\|^j_{m,j}+
\bigl\|\bigl(Y^{q,p},Z^{q,p}\bigr)\bigr\|^j_{m,j}
\\
&&\qquad\le C\bigl(\|\xi\|_{m+q,({(m+q-1)!}/{m!})j},\bigl(\bigl\|\partial^k_{\mathrm{sp}}f
\bigr\|_{\infty
}\bigr)_{k\le
m+q-1}\bigr),
\end{eqnarray*}
where $C$ is a constant depending on $\|\xi\|_{m+q,
({(m+q-1)!}/{m!})j}$ and
on\break $(\|\partial^k_{\mathrm{sp}}f\|_{\infty})_{k\le m+q-1}$.
\end{lemme}

\begin{prop}\label{prop1} Let $m \in\nset^{\star}$. Assume that
$\xi$
satisfies $\mathcal{H}_{m+q}^1$ and $f\in
C^{0,m+q-1,m+q-1}_b$. We recall $\mathcal{E}^{q,p}=\|
(Y^q-Y^{q,p},Z^q-Z^{q,p})\|^2_{\lp^2}$. We get
%
\begin{equation}
\label{eq7} \mathcal{E}^{q+1,p} \le C_1 T(T+1)
L_f^2 \mathcal{E}^{q,p} +\frac
{K_1(q,m)}{(p+1)\cdots(p+m)},
\end{equation}
where $C_1$ is a scalar and $K_1(q,m)$ depends on $T$, $m$, $\|\xi\|
_{m+q,2{(m+q-1)!}/{(m-1)!}}$ and
on $(\|\partial^k_{\mathrm{sp}} f\|_{\infty})_{1\le k \le m+q-1}$.


Since $\mathcal{E}^{0,p}=0$, we deduce from \eqref{eq7} that $
\mathcal{E}^{q,p}\le\frac{A_1(q,m)}{(p+1)^m}$ where $A_1(q,m):=
\frac{(C_1T(T+1)L_f^2)^q-1}{C_1T(T+1)L_f^2-1}K_1(q,m)$. Then,
$(Y^{p,q},Z^{p,q})$
converges to $(Y^q,Z^q)$ when $p$ tends to $\infty$ in
$\|(\cdot,\cdot)\|_{\lp^2}$; see \eqref{normeL2YZ} for the
definition of
the norm.
\end{prop}

\begin{rem} We deduce from Proposition~\ref{prop1} that for all $T$
and $L_f$,
we have $\lim_{p\rightarrow\infty}
\mathcal{E}^{q,p}=0$. When $C_1T(T+1) L_f^2<1$, that is, for $T$ small
enough, we
also get $\lim_{p\rightarrow\infty} \lim_{q\rightarrow\infty}
\mathcal{E}^{q,p}=0$.
\end{rem}

%
\begin{pf*}{Proof of Proposition~\ref{prop1}}
For the sake of
clearness, we
assume $d=1$. In the following, one notes $\Delta Y^{q,p}_t:=
Y^{q,p}_t-Y^{q}_t$, $\Delta Z^{q,p}_t:=
Z^{q,p}_t-Z^{q}_t$ and~$\Delta
f^{q,p}_t:=f(t,Y^{q,p}_t,Z^{q,p}_t)-f(t,Y^{q}_t,Z^{q}_t)$.
First,\vadjust{\goodbreak} we deal with\break  $\e[\sup_{0\le t \le T}
|\Delta Y^{q+1,p}_t|^2]$.
From \eqref{eqY} and \eqref{eqYpZp} we get
\begin{eqnarray*}
\Delta Y^{q+1,p}_t&=&\e_t\bigl[\cc_p
\bigl(F^{q,p}\bigr)-F^q\bigr]-\int_0^t
\Delta f^{q,p}_s \,ds
\\[-1pt]
&=&\e_t\bigl[\cc_p(\xi)-\xi\bigr]\\[-1pt]
&&{}+\e_t
\biggl[\cc_p \biggl(\int_0^T f
\bigl(s,Y^{q,p}_s,Z^{q,p}_s\bigr) \,ds
\biggr)-\int_0^T f\bigl(s,Y^{q}_s,Z^{q}_s
\bigr) \,ds \biggr]\\[-1pt]
&&{}-\int_0^t \Delta
f^{q,p}_s \,ds.
\end{eqnarray*}
We introduce $\pm\cc_p (\int_0^T f(s,Y^{q}_s,Z^{q}_s) \,ds
)$ in
the second conditional expectation. This leads to
\begin{eqnarray*}
\Delta Y^{q+1,p}_t&=&\e_t\bigl[\cc_p(
\xi)-\xi\bigr]+\e_t \biggl[\cc_p \biggl(\int
_0^T \Delta f^{q,p}_s \,ds
\biggr) \biggr]\\[-1pt]
&&{}+\e_t \biggl[\int_0^T
\cc_p \bigl(f\bigl(s,Y^{q}_s,Z^{q}_s
\bigr)\bigr) - f\bigl(s,Y^{q}_s,Z^{q}_s
\bigr) \,ds \biggr]
\\[-1pt]
&&{}-\int_0^t \Delta f^{q,p}_s
\,ds,
\end{eqnarray*}
where we have used the second property of Lemma~\ref{lem4} to
rewrite the third term.

From the previous equation, we bound $\e[\sup_{0\le t \le T} |\Delta
Y^{q+1,p}_t|^2]$ by
using Doob's inequality and the Lipschitz property of $f$
\begin{eqnarray*}
\e\Bigl[\sup_{0\le t \le T} \bigl| \Delta Y^{q+1,p}_t\bigr|^2
\Bigr]&\le&16 \e\bigl[\bigl|\cc_p(\xi)-\xi\bigr|^2\bigr]+16 \e \biggl[
\biggl\llvert \cc_p \biggl(\int_0^T
\Delta f^{q,p}_s \,ds \biggr)\biggr\rrvert ^2
\biggr]
\\[-1pt]
&&{}+16T\int_0^T \e \bigl[\bigl\llvert
\cc_p\bigl(f\bigl(s,Y^{q}_s,Z^{q}_s
\bigr)\bigr)- f\bigl(s,Y^{q}_s,Z^{q}_s
\bigr)\bigr\rrvert ^2 \bigr] \,ds\\[-1pt]
&&{}+8TL_f^2 \int
_0^T \e\bigl[\bigl|\Delta Y^{q,p}_s\bigr|^2+\bigl|
\Delta Z^{q,p}_s\bigr|^2\bigr] \,ds.
\end{eqnarray*}
To bound the second expectation of the previous inequality, we use the first
property of Lemma~\ref{lem4} and the Lispchitz property of $f$. Then we
bring together this term with the last one to get
%
\begin{eqnarray}
\label{eq5} \e\Bigl[\sup_{0\le t \le T} \bigl|\Delta Y^{q+1,p}_t\bigr|^2
\Bigr]&\le&16 \e\bigl[\bigl|\cc_p(\xi)-\xi\bigr|^2\bigr]\nonumber\\[-1pt]
&&{}+16T\int
_0^T \e \bigl[\bigl\llvert \cc_p
\bigl(f\bigl(s,Y^{q}_s,Z^{q}_s
\bigr)\bigr)- f\bigl(s,Y^{q}_s,Z^{q}_s
\bigr)\bigr\rrvert ^2 \bigr] \,ds
\\[-1pt]
&&{}+40TL_f^2 \int_0^T
\e\bigl[\bigl|\Delta Y^{q,p}_s\bigr|^2+\bigl|\Delta
Z^{q,p}_s\bigr|^2\bigr] \,ds.\nonumber\vadjust{\goodbreak}
\end{eqnarray}

Let us now upper bound $\e[\int_0^T |\Delta Z^{q+1,p}_s|^2 \,ds]$. To do
so, we use the It\^{o} isometry $\e[\int_0^T |\Delta Z^{q+1,p}_s|^2
\,ds]=\e[(\int_0^T
\Delta Z^{q+1,p}_s \,dB_s)^2]$. Using the definitions
\eqref{eqY}--\eqref{eqZp} of $Z^{q+1}_t$ and
$Z^{q+1,p}_t$ and the Clark--Ocone theorem leads to
\begin{eqnarray*}
\int_0^T \Delta Z^{q+1,p}_s
\,dB_s&=&F^q-\e\bigl(F^q\bigr)-\bigl(
\cc_p\bigl(F^{q,p}\bigr)-\e\bigl(\cc _p
\bigl(F^{q,p}\bigr)\bigr)\bigr)
\\
&=&Y^{q+1}_T+\int_0^T f
\bigl(s,Y_s^q,Z^q_s
\bigr)\,ds\\
&&{}-Y^{q+1}_0- \biggl(Y^{q+1,p}_T+
\int_0^T f\bigl(s,Y_s^{q,p},Z^{q,p}_s
\bigr)\,ds-Y^{q+1,p}_0 \biggr).
\end{eqnarray*}
Rearranging this summation makes
$\Delta Y^{q+1,p}_T - (\Delta Y^{q+1,p}_0)$ appear. We get
%
\begin{eqnarray}
\label{eq6} &&\e \biggl[\int_0^T \bigl|\Delta
Z^{q+1,p}_s\bigr|^2 \,ds \biggr]
\nonumber
\\
&&\qquad\le6\e\Bigl[\sup
_{0\le t
\le T}\bigl |\Delta Y^{q+1,p}_t\bigr|^2
\Bigr]\\
&&\qquad\quad{}+6TL_f^2 \int_0^T
\e\bigl[\bigl|\Delta Y^{q,p}_s\bigr|^2+\bigl|\Delta
Z^{q,p}_s\bigr|^2\bigr] \,ds.\nonumber
\end{eqnarray}

Since $\int_0^T
\e[|\Delta Y^{q,p}_s|^2+|\Delta Z^{q,p}_s|^2] \,ds \le(T+1)
\mathcal{E}^{q,p}$, by computing 7$\times$\eqref{eq5}${}+{}$\eqref{eq6}
we obtain
\begin{eqnarray*}
\mathcal{E}^{q+1,p} &\le&112 \e\bigl[\bigl|\cc_p(\xi)-
\xi\bigr|^2\bigr]\\
&&{}+112T\int_0^T \e \bigl[
\bigl\llvert \cc_p\bigl(f\bigl(s,Y^{q}_s,Z^{q}_s
\bigr)\bigr)- f\bigl(s,Y^{q}_s,Z^{q}_s
\bigr)\bigr\rrvert ^2 \bigr] \,ds\\
&&{}+286T(T+1)L_f^2
\mathcal{E}^{q,p}.
\end{eqnarray*}
Since $\xi$ and $f(s,Y^q_s,Z^q_s)$ belong to $\mathbb{D}^{m,2}$ ($\xi$
satisfies $\mathcal{H}^1_{m+q}$, $f \in C^{0,m+q-1,m+q-1}_b$ and
$(Y^q,Z^q) \in\mathcal{S}^{m,\infty}$ [see Lemma~\ref{lem5})], Lemma~\ref{lem2} gives
\begin{eqnarray*}
\mathcal{E}^{q+1,p} &\le& \frac{112}{(p+1)\cdots(p+m)}\bigl\|D^m \xi
\bigr\|^2_{\lp^2(\Omega
\times[0,T]^m)}
\\
&&{}+\frac{112T}{(p+1)\cdots(p+m)} \biggl(\int_0^T
\bigl\|D^m f\bigl(s,Y_s^q,Z^q_s
\bigr)\bigr\|^2_{\lp^2(\Omega
\times[0,T]^m)} \,ds \biggr)\\
&&{}+286T(T+1)L_f^2
\mathcal{E}^{q,p}.
\end{eqnarray*}
Since $\int_0^T\|D^m f(s,Y_s^q,Z^q_s)\|^2_{\lp^2(\Omega
\times[0,T]^m)} \,ds$ is bounded by $C(T,m, (\|\partial^k_{\mathrm{sp}}
f\|_{\infty})_{k \le m},\break\|(Y^q,Z^q)\|^{2m}_{m,2m})$, Lemma~\ref
{lem5} gives the result.
\end{pf*}

\subsection{Error due to the truncation of the basis}\label
{sectruncationbasis}

We are now interested in bounding the error between $(Y^{q,p},Z^{q,p})$
[defined by \eqref{eqYpZp}] and $(Y^{q,p,N},Z^{q,p,N})$ [defined by
\eqref{eqYZqpN}].

Before giving an upper bound for the error, we measure the error between
$\cc_p$ and $\cc_p^N$ for a r.v. satisfying \eqref{eq31} when $r= p$.

\begin{rem}\label{rem5}
Let $r\in\nset^{\star}$, $\xi$ satisfies $\mathcal{H}_{r+q}$ and $f
\in
C^{0,r+q-1,r+q-1}_b$. Then, for all integers $p$ and $q$,
$I_{q,p}:=\int_0^T
f(s,Y^{q,p}_s,Z^{q,p}_s) \,ds$ satisfies \eqref{eq31}; that is, for all
multi-index $\alpha$ such that $|\alpha|=r$, we have
\[
\bigl|\e\bigl(D^{\alpha}_{t_1,\ldots,t_r}I_{q,p}\bigr)-\e
\bigl(D^{\alpha}_{s_1,\ldots
,s_r}I_{q,p}\bigr)\bigr|\le
K^{I_{q,p}}_r\bigl((t_1-s_1)^{\beta_{I_{q,p}}}+
\cdots+(t_r-s_r)^{\beta
_{I_{q,p}}}\bigr),
\]
where $\beta_{I_{q,p}}=\frac{1}{2}\wedge\beta_{\xi}$ and
$K^{I_{q,p}}_r$ depends
on $K^{\xi}_r$, $\|\xi\|_{r+q,2{(r+q-1)!}/{(r-1)!}}$, $T$ and on
$(\|\partial^k_{\mathrm{sp}} f\|_{\infty})_{1\le k \le r+q-1}$.
\end{rem}




We refer to Section~\ref{proofrem5} for the proof of Remark~\ref{rem5}.
%
\begin{lemme}\label{lem6} Let $F$ denote a r.v. in $\lp^2(\m F_T)$ satisfying
\eqref{eq31} for $r=p$. We have
\[
\e\bigl(\bigl|\bigl(\cc_p^N -\cc_p\bigr)
(F)\bigr|^2\bigr) \le\bigl(K^{F}_p
\bigr)^2 \biggl(\frac{T}{N} \biggr)^{2\beta_F}\sum
_{i=1}^p i^2\frac{T^i}{i!}\le
\bigl(K^{F}_p\bigr)^2 \biggl(\frac{T}{N}
\biggr)^{2\beta_F} T(1+T)e^T,
\]
where $K^F_p$ and $\beta_F$ are defined in Hypothesis \ref{hypo3}.
\end{lemme}

We refer to Section~\ref{prooflem6} for the proof of the lemma.
%
\begin{prop}\label{prop4} Assume that $\xi$ satisfies $\mathcal
{H}_{p+q}$ and $f\in
C^{0,p+q-1,p+q-1}_b$. We recall
$\mathcal{E}^{q,p,N}:=\|(Y^{q,p}-Y^{q,p,N},Z^{q,p}-Z^{q,p,N})\|^2_{\lp^2}$.
We get
%
\begin{equation}
\label{eq8} \mathcal{E}^{q+1,p,N} \le C_2 T(T+1)
L_f^2 \mathcal{E}^{q,p,N} +K_2(q,p)
\biggl(\frac{T}{N} \biggr)^{1 \wedge2\beta_\xi},
\end{equation}
where $C_2$ is a scalar and $K_2(q,p)$ depends
on $K^{\xi}_p$, $T$, $\|\xi\|_{p+q,2{(p+q-1)!}/{(p-1)!}}$ and on
$(\|\partial^k_{\mathrm{sp}}
f\|_{\infty})_{1\le k \le p+q-1}$.

Since $\mathcal{E}^{0,p,N}=0$, we deduce from \eqref{eq8} that $
\mathcal{E}^{q,p,N}\le A_2(q,p)  (\frac{T}{N}  )^{1 \wedge
2\beta_\xi}$, where $A_2(q,p):=
K_2(q,p)T(T+1)e^T\frac{(C_2T(T+1)L_f^2)^q-1}{C_2T(T+1)L_f^2-1}$. Then,
$(Y^{p,q,N},Z^{p,q,N})$
converges to $(Y^{q,p},Z^{q,p})$ when $N$ tends to $\infty$ in
$\|(\cdot,\cdot)\|_{\lp^2}$.
\end{prop}



\begin{pf*}{Proof of Proposition~\ref{prop4}}
For the sake of clarity, we
assume $d=1$. In the following, we note $\Delta Y^{q,p,N}_t:=
Y^{q,p,N}_t-Y^{q,p}_t$, $\Delta Z^{q,p,N}_t:=
Z^{q,p,N}_t-Z^{q,p}_t$ and $\Delta
f^{q,p,N}_t:=f(t,Y^{q,p,N}_t,Z^{q,p,N}_t)-f(t,Y^{q,p}_t,Z^{q,p}_t)$.
First, we deal with $\e[\sup_{0\le t \le T}
|\Delta Y^{q+1,p,N}_t|^2]$.
From \eqref{eqYpZp} and \eqref{eqYZqpN} we get
\[
\Delta Y^{q+1,p,N}_t=\e_t\bigl[
\cc_p^N\bigl(F^{q,p,N}\bigr)-\cc_p
\bigl(F^{q,p}\bigr)\bigr]+\int_0^t
\Delta f^{q,p,N}_s \,ds.
\]
Following the same steps as in the proof of Proposition~\ref{prop1},
we get
%
\begin{eqnarray}
\label{eq10} &&\e\Bigl[\sup_{0\le t \le T} \bigl|\Delta Y^{q+1,p,N}_t\bigr|^2
\Bigr]\nonumber\\
&&\qquad\le16 \e\bigl[\bigl|\cc_p^N(\xi)-\cc_p(
\xi)\bigr|^2\bigr]
\nonumber
\\[-8pt]
\\[-8pt]
\nonumber
&&\qquad\quad{}+16\e \biggl[\biggl\llvert \bigl(\cc_p^N
- \cc_p\bigr) \biggl(\int_0^T f
\bigl(s,Y^{q,p}_s,Z^{q,p}_s\bigr)\,ds
\biggr)\biggr\rrvert ^2 \biggr]
\\
&&\qquad\quad{}+40TL_f^2 \int_0^T
\e\bigl[\bigl|\Delta Y^{q,p,N}_s\bigr|^2+\bigl|\Delta
Z^{q,p,N}_s\bigr|^2\bigr] \,ds.\nonumber
\end{eqnarray}

Let us now upper bound $\e[\int_0^T |\Delta Z^{q+1,p,N}_s|^2
\,ds]$. Following the same steps as in the proof of Proposition~\ref{prop1},
we get
%
\begin{eqnarray}
\label{eq9}&& \e \biggl[\int_0^T \bigl|\Delta
Z^{q+1,p,N}_s\bigr|^2 \,ds \biggr]\nonumber\\
&&\qquad\le 6\e\Bigl[\sup
_{0\le t
\le T} \biggl|\Delta Y^{q+1,p,N}_t\bigr|^2
\Bigr]
\\
&&\qquad\quad{}+6TL_f^2 \int_0^T
\e\bigl[\bigl|\Delta Y^{q,p,N}_s\bigr|^2+\bigl|\Delta
Z^{q,p,N}_s\bigr|^2\bigr] \,ds.\nonumber
\end{eqnarray}

Adding $7\times$\eqref{eq10} and \eqref{eq9} gives
\begin{eqnarray*}
\mathcal{E}^{q+1,p,N}& \le& 112 \e\bigl[\bigl|\bigl(\cc_p^N
- \cc_p\bigr) (\xi)\bigr|^2\bigr]\\
&&{}+112\e \biggl[\biggl\llvert
\bigl(\cc_p^N - \cc_p\bigr) \biggl(\int
_0^T f\bigl(s,Y^{q,p}_s,Z^{q,p}_s
\bigr)\,ds \biggr)\biggr\rrvert ^2 \biggr]
\\
&&{}+286T(T+1)L_f^2 \mathcal{E}^{q,p,N}.
\end{eqnarray*}

Since $\xi$ and $I_{q,p}$ satisfy \eqref{eq31} (see Remarks \ref{rem9}
and \ref{rem5}), Lemma~\ref{lem6} gives
\begin{eqnarray*}
\mathcal{E}^{q+1,p,N}& \le&112 \biggl(\frac{T}{N} \biggr)^{2\alpha_\xi
\wedge1}
T (T+1) e^T \bigl(\bigl(K^{\xi}_p
\bigr)^2+\bigl(K^{I_{q,p}}_p\bigr)^2
\bigr)\\
&&{}+286T(T+1)L_f^2 \mathcal{E}^{q,p,N},
\end{eqnarray*}
and \eqref{eq8} follows.
\end{pf*}

\subsection{Error due to the Monte Carlo approximation}\label{sectMC}

We are now interested in bounding the error between $(Y^{q,p,N},Z^{q,p,N})$
defined by \eqref{eqYZqpN} and $(Y^{q,p,N,M},\allowbreak Z^{q,p,N,M})$ defined by
\eqref{eqYZqpNM}. $\cc_p^{N,M}$ is defined by \eqref{dchapeau} and
\eqref{chaosdecMC}. In this section, we assume that the coefficients
$\hat{d^n_k}$ are independent of the vector $(G_1,\ldots,G_N)$, which
corresponds to the second approach proposed in Remark~\ref{rem2}.

Before giving an upper bound for the error, we measure the error between
$\cc_p^N$ and $\cc_p^{N,M}$ for a r.v. satisfying $\mathcal
{H}^3_{p,N}$ (see
Hypothesis \ref{hypo4}).


\begin{lemme}\label{lem7}
Let $F$ be a r.v. satisfying Hypothesis $\mathcal{H}^3_{p,N}$. We have
\[
\e\bigl(\bigl|\bigl(\cc_p^N-\cc_p^{N,M}
\bigr) (F)\bigr|^2\bigr) = \frac{1}{M}V_{p,N}(F).
\]
Moreover, we have $ \e(|\cc_p^{N,M}(F)|^2) \le\e(|F|^2)+\frac
{1}{M}V_{p,N}(F)$.
\end{lemme}
We refer to Section~\ref{prooflem7} for the proof of the lemma.

\begin{prop}\label{prop5} Let $\xi$ satisfy Hypothesis $\mathcal
{H}^3_{p,N}$ and $f$ be a bounded function. Let $\mathcal
{E}^{q,p,N,M}:=\|(Y^{q,p,N}-Y^{q,p,N,M},Z^{q,p,N}-Z^{q,p,N,M})\|^2_{\lp^2}$.
We get
\[
\mathcal{E}^{q+1,p,N,M} \le C_3 T(T+1) L_f^2
\mathcal{E}^{q,p,N,M} +\frac{K_3(q,p,N)}{M},
\]
where $C_3$ is a scalar and $K_3(q,p,N):=168 (V_{p,N}(\xi)+T^2
\|f\|^2_{\infty} \sum_{k=0}^p {N \choose k}  )$.

Since $\mathcal{E}^{0,p,N,M}=0$, we deduce from the previous
inequality that $
\mathcal{E}^{q,p,N,M}\le\frac{A_3(q,p,N)}{M} $, where $A_3(q,p,N):=
K_3(q,p,N)\frac{(C_3T(T+1)L_f^2)^q-1}{C_3T(T+1)L_f^2-1}$. Then
$(Y^{p,q,N,M},\allowbreak Z^{p,q,N,M})$
converges to $(Y^{q,p,N},Z^{q,p,N})$ when $M$ tends to $\infty$ in
$\|(\cdot,\cdot)\|_{\lp^2}$.
\end{prop}

\begin{pf*}{Proof of Proposition~\ref{prop5}}
For the sake of clarity, we
assume $d=1$. In the following, note that $\Delta Y^{q,p,N,M}_t:=
{Y}^{q,p,N,M}_t-Y^{q,p,N}_t$, $\Delta Z^{q,p,N,M}_t:=
{Z}^{q,p,N,M}_t-Z^{q,p,N}_t$ and $\Delta
f^{q,p,N,M}_t:=f(t,{Y}^{q,p,N,M}_t,{Z}^{q,p,N,M}_t)-f(t,Y^{q,p,N}_t,\break Z^{q,p,N}_t)$.
First, we deal with $\e[\sup_{0\le t \le T}
|\Delta Y^{q+1,p,N,M}_t|^2]$.
From \eqref{eqYZqpN} and \eqref{eqYZqpNM} we get
\[
\Delta Y^{q+1,p,N,M}_t=\e_t\bigl[
\cc_p^{N,M}\bigl({F}^{q,p,N,M}\bigr)-\cc
^N_p\bigl(F^{q,p,N}\bigr)\bigr]+\int
_0^t \Delta f^{q,p,N,M}_s \,ds.
\]
By introducing $\pm\cc_p^{N}({F}^{q,p,N,M})$ and by using Lemma~\ref{lem9}, we obtain
\begin{eqnarray*}
\e\Bigl[\sup_{0\le t \le T}\bigl |\Delta Y^{q+1,p,N,M}_t\bigr|^2
\Bigr]&\le&12 \e\bigl[\bigl|\bigl(\cc_p^{N,M}-\cc^N_p
\bigr) \bigl({F}^{q,p,N,M}\bigr)\bigr|^2\bigr]\\
&&{}+12\e
\bigl(\bigl|{F}^{q,p,N,M}-F^{q,p,N}\bigr|^2 \bigr)
\\
&&{}+6TL_f^2 \int_0^T
\e\bigl[\bigl|\Delta Y^{q,p,N,M}_s\bigr|^2+\bigl|\Delta
Z^{q,p,N,M}_s\bigr|^2\bigr] \,ds.
\end{eqnarray*}
From Lemma~\ref{lem7}, we get
$\e[|(\cc_p^{N,M}-\cc^N_p)({F}^{q,p,N,M})|^2] \le
\frac{2}{M} (V_{p,N}(\xi)+\break V_{p,N} (\int_0^T
f({\theta}^{q,p,N,M}_s) \,ds  ) )$. Then
%
\begin{eqnarray}
\label{eq22} &&\e\Bigl[\sup_{0\le t \le T} \bigl|\Delta Y^{q+1,p,N,M}_t\bigr|^2
\Bigr]\nonumber\\
&&\qquad\le \frac{24}{M} \Biggl(V_{p,N}(\xi)+T^2 \|f
\|^2_{\infty} \sum_{k=0}^p
\pmatrix{N \cr k} \Biggr)
\\
&&\qquad\quad{}+30TL_f^2 \int_0^T
\e\bigl[\bigl|\Delta Y^{q,p,N,M}_s\bigr|^2+\bigl|\Delta
Z^{q,p,N,M}_s\bigr|^2\bigr] \,ds.\nonumber
\end{eqnarray}

Let us now upper bound $\e[\int_0^T |\Delta Z^{q+1,p,N,M}_s|^2
\,ds]$. Following the same steps as in the proof of Proposition~\ref{prop1},
we get
%
\begin{eqnarray}
\label{eq21} &&\e \biggl[\int_0^T\bigl |\Delta
Z^{q+1,p,N,M}_s\bigr|^2 \,ds \biggr]\nonumber\\
&&\qquad\le 6\e \Bigl[\sup
_{0\le t
\le T}\bigl |\Delta Y^{q+1,p,N,M}_t\bigr|^2
\Bigr]\\
&&\qquad\quad{}+6TL_f^2 \int_0^T
\e\bigl[\bigl|\Delta Y^{q,p,N,M}_s\bigr|^2+\bigl|\Delta
Z^{q,p,N,M}_s\bigr|^2\bigr] \,ds.\nonumber
\end{eqnarray}

Adding $7\times$\eqref{eq22} and \eqref{eq21} gives the result.
\end{pf*}

\section{Numerical examples}
\label{secnumex}
The computations have been done on a PC INTEL Core 2 Duo P9600 2.53 GHz with
4Gb of RAM.

\subsection{Nonlinear driver and path-dependent terminal condition}
We consider the case $d=1$, $f(t,y,z)=\cos(y)$ and $\xi=\sup_{0\le t
\le1} B_t$.
\begin{itemize}
\item\textbf{Convergence in $p$.}
Tables~\ref{tab1} and \ref{tab2} represent the evolution of $\overline
{Y}^{q,p,N,M}_0$ and
$\overline{Z}^{q,p,N,M}_0$ w.r.t $q$ (Picard's
iteration index), when
$p=2$ and $p=3$. We also give the CPU time needed to get $\overline
{Y}^{6,p,N,M}_0$ and
$\overline{Z}^{6,p,N,M}_0$. We fix $M=10^5$ and $N=20$. The seed of
the generator is
also fixed.
%
\begin{table}[b]
\tabcolsep=0pt
\caption{Evolution of $\overline{Y}^{q,p,N,M}_0$ w.r.t. Picard's
iterations, $M=10^5$, $N=20$ and
CPU time}
\label{tab1}
\begin{tabular*}{\textwidth}{@{\extracolsep{\fill}}lccccccc@{}}
\hline
\textbf{Iterations} & \textbf{1} & \textbf{2} & \textbf{3} &
\textbf{4} & \textbf{5} & \textbf{6} & \textbf{CPU time}\\
\hline
$p=2$ & 1.656357 & 1.017117 & 1.237135 & 1.186691 & 1.195462 &
1.194256 & \phantom{0}14.06\\
$p=3$ & 1.656357 & 1.012091 & 1.234398 & 1.183544 & 1.192367 &
1.191173 & 174.09\\
\hline
\end{tabular*}
\end{table}

\begin{table}
\tabcolsep=0pt
\caption{Evolution of $\overline{Z}^{q,p,N,M}_0$ w.r.t. Picard's
iterations, $M=10^5$, $N=20$ and
CPU time}
\label{tab2}
\begin{tabular*}{\textwidth}{@{\extracolsep{\fill}}lccccccc@{}}
\hline
\textbf{Iterations} & \textbf{1} & \textbf{2} & \textbf{3} &
\textbf{4} & \textbf{5} & \textbf{6} & \textbf{CPU time}\\
\hline
$p=2$ & 0.969128 & 0.249148 & 0.525273 & 0.459326 & 0.470069 &
0.469117 & \phantom{0}14.06\\
$p=3$ & 0.969128 & 0.242977 & 0.523846 & 0.455827 & 0.466903 &
0.465939 & 174.09\\
\hline
\end{tabular*}
\end{table}

Note that the difference between the values of $\overline
{Y}^{q,2,N,M}_0$ and
$\overline{Y}^{q,3,N,M}_0$ (resp., $\overline{Z}^{q,2,N,M}_0$ and
$\overline{Z}^{q,3,N,M}_0$) does not exceed $0.2 \%$ (resp., $0.6 \%$).
This is due to the fast convergence of the
algorithm in $p$. The CPU time is $12$ times higher when $p=3$ than
when $p=2$. Then, the use of order $3$ in the chaos decomposition is not
necessary. In the following, we take $p=2$.
\item\textbf{Convergence in $M$.}
Figure~\ref{fig1} illustrates the evolution of $\overline
{Y}^{q,p,N,M}_0$ and $\overline{Z}^{q,p,N,M}_0$ w.r.t. $q$ when
$p=2$ and $N=20$ for different values of $M$. The seed of the generator
is random. When $M$ equals $10^4$ and $10^5$ the algorithm stabilizes
after very few iterations. When
$M=10^3$, there is no convergence.

\begin{figure}[b]

\includegraphics{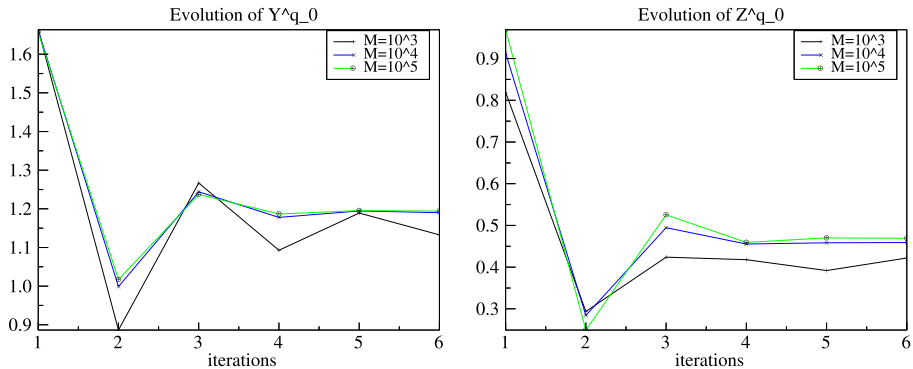}

\caption{Evolution of $\overline{Y}^{q,p,N,M}_0$ and $\overline
{Z}^{q,p,N,M}_0$ w.r.t. $q$ and $M$
when $N=20$, $p=2 -
\xi=\sup_{0\le t \le1} B_t$, $f(t,y,z)=\cos(y)$.}
\label{fig1}
\end{figure}

\begin{figure}

\includegraphics{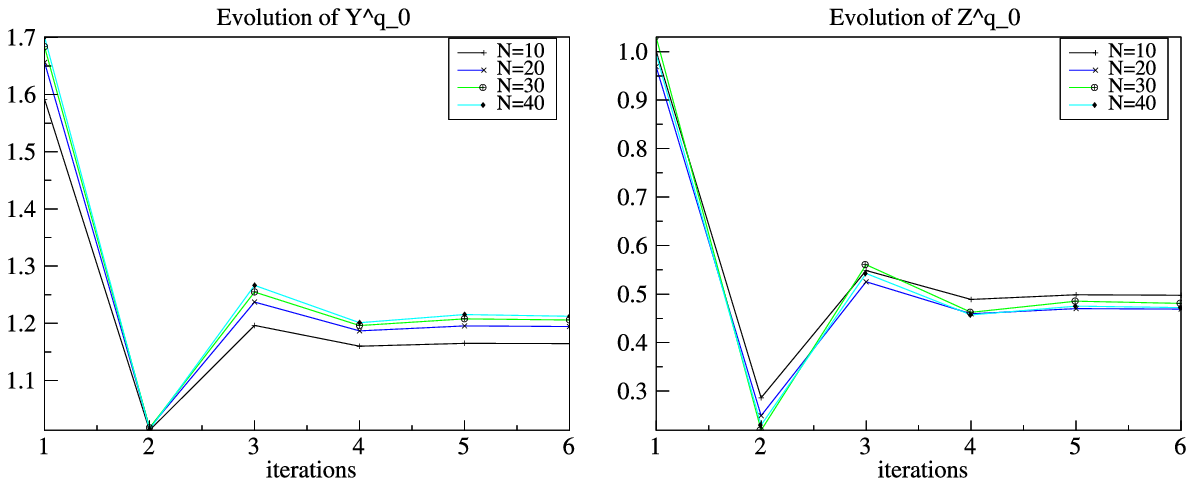}

\caption{Evolution of $\overline{Y}^{q,p,N,M}_0$ and $\overline
{Z}^{q,p,N,M}_0$ w.r.t. $N$ when $M=10^5$, $p=2 -
\xi=\sup_{0\le t \le1} B_t$, $f(t,y,z)=\cos(y)$.}
\label{fig2}
\end{figure}

\item\textbf{Convergence in $N$.}
Figure~\ref{fig2} illustrates the evolution of $\overline
{Y}^{q,p,N,M}_0$ and $\overline{Z}^{q,p,N,M}_0$ w.r.t. $q$ when
$p=2$ and $M=10^5$ for different values of $N$. The seed of the generator
is random. The algorithm converges even when $N=10$, but $\overline
{Y}^{6,p,10,M}_0$ is quite
below $\overline{Y}^{6,p,40,M}_0$.
\end{itemize}

\subsection{Linear driver-financial benchmark}
We consider the case of pricing and hedging a discrete down and out barrier
call option, that is, $f(t,y,z)=-ry$ and $\xi:=(S_T-K)_+\ind_{\forall n
\in
[0,N] S_{t_n} \ge L}$, where $S$ represents the Black--Scholes diffusion
\[
S_t=S_0e^{(r-({1}/{2})\sigma^2)t+\sigma W_t} \qquad\forall t\in[0,T].
\]

The option parameters are $r=0.01$, $\sigma=0.2$, $T=1$, $K=0.9$, $L=0.85$,
$S_0=1$ and $N=20$ ($N$ is also the number of time discretizations of
the chaos decomposition).

\begin{figure}[b]

\includegraphics{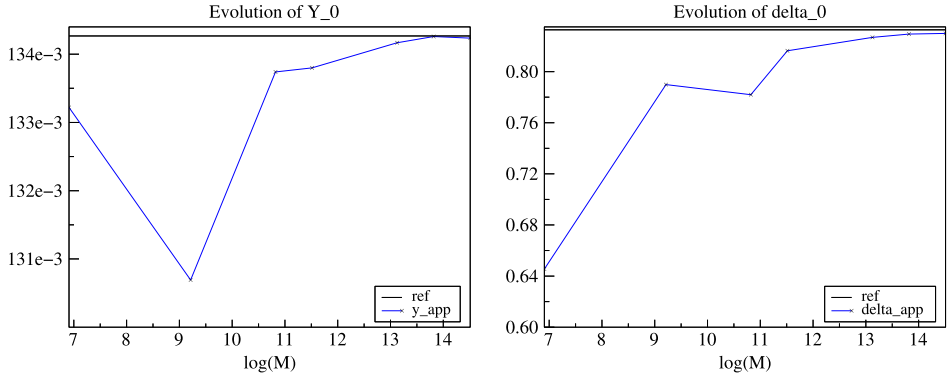}

\caption{Evolution of $\overline{Y}^{q,p,N,M}_0$ and
$\delta_0:=\frac{\overline{Z}^{q,p,N,M}_0}{\sigma S_0}$ w.r.t.
$log(M)$ when
$N=20$, $p=2$, $q=5$, discrete down and out barrier call option.}
\label{fig3}
\end{figure}
We can get a benchmark for $Y_0$ and $Z_0$ by using a variance
reduction Monte Carlo method. For this set of parameters, the reference values
are $Y_0=0.134267$ with a confidence interval $7.9468e-05$ and
$\delta_0=\frac{Z_0}{\sigma S_0}=0.8327$. We compare them with
$\overline{Y}^{q,p,N,M}_0$
and $\frac{\overline{Z}^{q,p,N,M}_0}{\sigma S_0}$ when $N=20$, $p=2$,
$q=5$ (we choose
the first value of $q$ from which the algorithm has converged) for different
values of $M$. Figure~\ref{fig3} represents the evolution of
$\overline
{Y}^{5,p,N,M}_0$
and $\delta^{5,p,N,M}_0$ w.r.t. $\log(M)$. Notice that for $M=10^6$ the
computed values are very close to the reference values.

\begin{figure}

\includegraphics{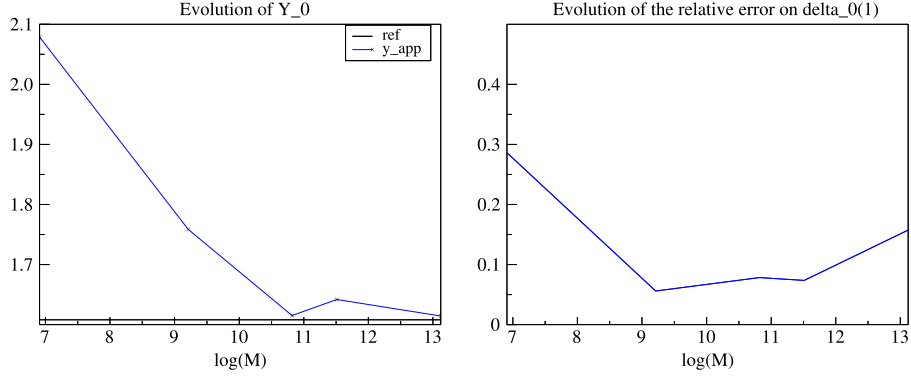}

\caption{Evolution of $\overline{Y}^{q,p,N,M}_0$ and $\delta_0(1)$
w.r.t. $\log(M)$
when $N=20$, $p=2$, $q=5$, $d=5$, basket put option with different
interest and borrowing rates.}
\label{fig4}
\end{figure}

\subsection{Nonlinear driver in dimension $5$, financial benchmark}
We consider the pricing and hedging of a put basket option in dimension $5$,
that is, $\xi=(K-\frac{1}{5}\sum_{i=1}^5 S_T^i)_+$, where
\[
\forall i=1,\ldots,5\qquad  S^i_t=S^i_0
e^{(\mu^i-{(\sigma
^i)^2}/{2})t+\sigma^i B^i_t}.
\]
$\mu^i$ (resp., $\sigma^i$) represents the trend (resp., the
volatility) of the $i$th asset. $B=(B^1,\ldots,B^5)$ is a $5$-dimensional
Brownian motion such that $\langle B^i, B^j \rangle_t = \rho t \ind_{i
\neq j}
+ t \ind_{i = j}$. We suppose that $\rho\in(-\frac{1}{4}, 1)$,
which ensures
that the matrix $C=(\rho\ind_{i \neq j} +\ind_{i = j})_{1\leq
i,j\leq
5}$ is
positive definite. We also assume that the borrowing rate $R$ is higher than
the bond rate $r$. In such a case, pricing and hedging the put basket
option is
equivalent to solving a BSDE with terminal condition $\xi$ and with
driver $f$
defined by $f(t,y,z)=-ry- \theta\cdot z +(R-r)(y-\sum_{i=1}^5
(\Sigma^{-1}z)_i)^-$, where $\theta:=\Sigma^{-1}(\mu-r{\mathbf1})$
(${\mathbf1}$ is the vector whose every component is one), and $\Sigma
$ is the
matrix defined by $\Sigma_{ij}=\sigma^i L_{ij}$ ($L$ denote the lower
triangular matrix involved in the Cholesky decomposition $C=LL^*$). We refer
to \cite{karoui97}, Example~1.1, for more details.

The option parameters are $r=0.02$, $R=0.1$, $T=1$, $K=95$, $\rho
=0.1$, and
for all $i=1,\ldots,5$, $S_0^i=100$, $\mu_0^i=0.05$ and $\sigma
_0^i=0.2$. Figure~\ref{fig4}
represents the evolution of $\overline{Y}^{5,p,N,M}_0$, the
approximated price at time
$0$ and the relative error on $\delta_0^1:=\frac{(\Sigma
^{-1}\overline
{Z}^{5,p,N,M}_0)^1}{S_0^1}$---the quantity of asset $1$ to possess at time
$0$---w.r.t. $\log(M)$. We compare our results with the ones obtained using the
Algorithm proposed in~\cite{GL10} (cited here as reference values).
The CPU time needed to compute price and delta when $M=50\mbox{,}000$ and
$N=20$ is $161$~s. Notice that the convergence is very fast and quite
accurate for $M=50\mbox{,}000$.

\textit{Conclusion}.
In this paper, we use Wiener chaos expansions together with the Picard
procedure to compute the solution to \eqref{eqmain}. Once the chaos
decomposition of $F^q$ is computed, we get explicit formulas
for both conditional expectations and the Malliavin derivative of conditional
expectations. This enables us to easily compute $(Y^q,Z^q)$.
Numerically, we
obtain fast and accurate results, which encourage us to extend these
results to
other type of BSDEs, like 2-BSDEs. It is also possible to
couple these Wiener chaos expansions together with the dynamic programming
approach. This will be the subject of a forthcoming publication.

\begin{appendix}\label{app}
\section{Technical results of Section \lowercase{\protect\texorpdfstring{\ref{secconvres}}{4}}}\label
{secappendix}

In the following, for any regular r.v. $F \in\mathcal{F}_T$,
$D^{(l_0)}_t\Delta_iD^{(l_1)}_s F$ denotes
$D^{(l_0)}_{t_1,\ldots,t_{l_0}}\times\allowbreak (D^{(l_1+1)}_{t_i,s_{i+1},\ldots,s_{i+l_1}}
F-D^{(l_1+1)}_{s_i,\ldots,s_{i+l_1}}
F)$.
\subsection{Proof of Remark \protect\texorpdfstring{\ref{rem6}}{4.5}}\label{proofrem6}
Before proving Remark~\ref{rem6}, we prove the following lemma.
%
\begin{lemme}\label{lem18}
Let $X$ be the $\rset^n$-valued process solution of
\[
X_t=x+\int_0^t
b(s,X_s)\,ds +\int_0^t
\sigma(s,X_s)\,dB_s,
\]
where $B$ is a $d$-dimensional Brownian motion and $b\dvtx [0,T] \times
\rset^n
\rightarrow\rset^n$ and $\sigma\dvtx  [0,T] \times\rset^n \rightarrow
\rset^{n\times d}$ are two $C^{0,m}$ functions uniformly lipschitz
w.r.t. $x$ and H\"{o}lder continuous of parameter $\frac{1}{2}$ w.r.t.
$t$, with
linear growth in $x$ (of constant $K$) and with bounded derivatives. Then:
\begin{itemize}
\item$\forall l
\le m$, $\forall j\ge2$ we have
%
\begin{equation}
\label{eq23} M^j_l:=\sup_{t_1\le\cdots\le t_l} \e
\Bigl(\sup_{r \in[t_l, T]}\bigl |D^{(l)}_{t_1,\ldots,t_l}
X_r\bigr|^j\Bigr) < \infty,
\end{equation}
the upper bound depends on $(\|b^{(l')}\|_{\infty})_{l' \le l}$,
$(\|\sigma^{(l')}\|_{\infty})_{l' \le l}$, $x$ and $K$,
\item$\forall j \ge2$, $\forall i \in
\{1,\ldots,m\}, \forall l_0\le i-1$, $\forall l_1 \le m-i$, we have
%
\begin{equation}
\label{eq24} \sup_{t_1\le\cdots\le t_{l_0}}\sup_{s_{i+1}\le\cdots\le s_{i+l_1}} \e
\Bigl(\sup_{r\in[s_{i+l_1},T]}\bigl|D^{(l_0)}_t
\Delta_iD^{(l_1)}_s X_r
\bigr|^j\Bigr) < k^X_l(j) (t_i-s_i)^{{j}/{2}},\hspace*{-30pt}
\end{equation}
where $l:=l_0+l_1+1$ and $k^X_l$ depends on $T$, $(M_{l'}^{j'})_{l'\le
l,j' \le lj}$, $(\|b^{(l')}\|_{\infty})_{l' \le l}$, and on
$(\|\sigma^{(l')}\|_{\infty})_{l' \le l}$.
\end{itemize}
\end{lemme}

\begin{pf*}{Proof of Lemma~\ref{lem18}}
The first point is proved in
\cite{nualart06}, Theorem~2.2.2. For the sake of clarity, we prove
the second result for $d=1$. We also assume that the vectors
$(t_1,\ldots,t_n)$ and $(s_1,\ldots,s_n)$ are such that $0 \le s_1
\le t_1
\le s_2 \le\cdots\le s_n \le t_n \le T$. We do it by induction on
$l_0$ and
$l_1$. We detail the case $b$ and $\sigma$ only
depending on $x$ and do the proof for $l_0=l_1=0$ and $l_0=0$,
$l_1=1$. We recall that under these hypotheses on $b$ and $\sigma$, we have
$\forall l \le m$
\[
\sup_{t_1\le\cdots\le
t_{l}}\e\bigl[\bigl|D^{(l)}_{t_1,\ldots,t_{l}}(X_{t_{l+1}}-X_{s_{l+1}})\bigr|^j
\bigr] \le C(t_{l+1}-s_{l+1})^{{j}/{2}},
\]
where $C$ depends on $T$, $j$, $(M_{l'}^{j'})_{l'\le l,j' \le lj}$ and
on $(\|b^{(j')}\|_{\infty})_{j' \le j}$, and on
$(\|\sigma^{(j')}\|_{\infty})_{j' \le j}$.

\textit{Case} $l_0=l_1=0$. We have
\begin{eqnarray*}
D_{t_n}X_r&=&\int_{t_n}^r
b'(X_u)D_{t_n}X_u \,du +
\sigma(X_{t_n})\\
&&{}+\int_{t_n}^r
\sigma'(X_u)D_{t_n}X_u \,dBu.
\end{eqnarray*}
Then
\begin{eqnarray*}
\Delta_n X_r&:=& D_{t_n}X_r-D_{s_n}X_r
\\
&=&\int_{t_n}^r b'(X_u)
\Delta_n X_u \,du-\int_{s_n}^{t_n}
b(X_u) D_{s_n}(X_u)\,du
\\
&&{}+\sigma(X_{t_n})-\sigma(X_{s_n})+\int_{t_n}^r
\sigma'(X_u)\Delta_n X_u \,dBu\\
&&{}-
\int_{s_n}^{t_n} \sigma'(X_u)
D_{s_n}(X_u)\,dBu.
\end{eqnarray*}
In the following, $C$ denotes a generic constant depending only on $T$
and $j$, and
$L_{\sigma}$ denotes the Lipschitz contant of $\sigma$.
\begin{eqnarray*}
|\Delta_n X_r|^j &\le&C \biggl(
\bigl\|b'\bigr\|_{\infty}^j \int_{t_n}^r|
\Delta_n X_u|^j \,du+(t_n-s_n)^{j-1}
\bigl\|b'\bigr\|_{\infty}^j \int_{s_n}^{t_n}
\bigl| D_{s_n}(X_u)\bigr|^j\,du
\\
&&\hspace*{15pt}{}+L_{\sigma}^j |X_{t_n}-X_{s_n}|^j+
\biggl\llvert \int_{t_n}^r \sigma
'(X_u)\Delta_n X_u\, dBu\biggr
\rrvert ^j\\
&&\hspace*{156pt}{}+\biggl\llvert \int_{s_n}^{t_n}
\sigma'(X_u) D_{s_n}(X_u)\,dBu
\biggr\rrvert ^j \biggr).
\end{eqnarray*}
We introduce $\Psi^{0,j}_n(T):= \e[\sup_{r\in[t_n, T]} |\Delta_n
X_r|^j]$. Doob's inequality and the Burkh\"{o}lder--Davis--Gundy
inequality lead to
\begin{eqnarray*}
&&\Psi^{0,j}_n(T) \le C \biggl( \bigl(\bigl\|b'
\bigr\|_{\infty}^j+ \bigl\|\sigma'\bigr\|_{\infty}^j
\bigr) \int_{t_n}^T\Psi^{0,j}_n(u)
\,du+\bigl\|b'\bigr\|_{\infty}^jM_1^j(t_n-s_n)^{j}\\
&&\hspace*{170pt}{}+
\bigl(L_{\sigma}^j+ \bigl\|\sigma'\bigr\|
_{\infty}^jM_1^j\bigr)|t_n-s_n|^{{j}/{2}}
\biggr).
\end{eqnarray*}
Gronwall's lemma yields the result.\vadjust{\goodbreak}

\textit{Case} $l_0=0, l_1=1$. We consider $\Delta_{n-1}D_{t_n}X_r=
D_{t_{n-1},t_n}X_r - D_{s_{n-1},t_n}X_r$. We have
\begin{eqnarray*}
D_{t_{n-1},t_n}X_r&=&\int_{t_n}^r
b''(X_u)D_{t_{n-1}}X_u
D_{t_n}X_u + b'(X_u)
D_{t_{n-1},t_n}X_u \,du \\
&&{}+\sigma'(X_{t_n})D_{t_{n-1}}X_{t_n}
\\
&&{}+\int_{t_n}^r \sigma''(X_u)D_{t_{n-1}}X_u
D_{t_n}X_u + \sigma'(X_u)
D_{t_{n-1},t_n}X_u \,dBu.
\end{eqnarray*}
Then
\begin{eqnarray*}
\Delta_{n-1}D_{t_n}X_r&=&\int
_{t_n}^r b''(X_u)
\Delta_{n-1}X_u D_{t_n}X_u +
b'(X_u) \Delta_{n-1}D_{t_n}X_u
\,du \\
&&{}+\sigma'(X_{t_n})\Delta _{n-1}X_{t_n}
\\
&&{}+\int_{t_n}^r \sigma''(X_u)
\Delta_{n-1}X_u D_{t_n}X_u +
\sigma'(X_u) \Delta_{n-1}D_{t_n}X_u\,dBu.
\end{eqnarray*}
Doob's inequality and the Burkh\"{o}lder--Davis--Gundy inequality lead to
\begin{eqnarray*}
\hspace*{-4pt}&&\e\Bigl[\sup_{r \in[t_n,T]}\bigl|\Delta_{n-1}D_{t_n}X_r\bigr|^j
\Bigr]\\
\hspace*{-4pt}&&\qquad\le C \biggl( \int_{t_n}^T
\bigl\|b''\bigr\|_{\infty}^j \e\bigl[\bigl|
\Delta_{n-1}X_u\bigr|^j |D_{t_n}X_u|^j
\bigr] \\
\hspace*{-4pt}&&\hspace*{65pt}{}+ \bigl\|b'\bigr\|_{\infty}^j \e\bigl[|
\Delta_{n-1}D_{t_n}X_u|^j\bigr] \,du
\\
\hspace*{-4pt}&&\hspace*{47pt}{}+\bigl\|\sigma'\bigr\|^j_{\infty}\e\bigl[|
\Delta_{n-1}X_{t_n}|^j\bigr]\\
\hspace*{-4pt}&&\hspace*{47pt}{}+\int
_{t_n}^T \bigl\|\sigma''
\bigr\|^j_{\infty}\e\bigl[|\Delta_{n-1}X_u|^j|
D_{t_n}X_u|^j\bigr] + \bigl\| \sigma
'\bigr\|^j_{\infty}\e\bigl[|\Delta_{n-1}D_{t_n}X_u|^j
\bigr] \,du \biggr).
\end{eqnarray*}

We introduce $\Psi^{1,j}_{n-1}(T):=\sup_{t_n \le T} \e[\sup_{r\in
[t_n,T]} |\Delta_{n-1}D_{t_n}
X_r|^j]$. The Cau\-chy--Schwarz inequality yields
\begin{eqnarray*}
\Psi^{1,j}_{n-1}(T)& \le& C \biggl(\bigl(\bigl\|b'
\bigr\|_{\infty}^j+\bigl\|\sigma'\bigr\| _{\infty
}^j
\bigr)\int_{t_n}^T \Psi^{1,j}_{n-1}(u)\,du\\
&&\hspace*{15pt}{}+
\bigl(\bigl\|b''\bigr\|_{\infty}^j+\bigl\|
\sigma''\bigr\|^j_{\infty}\bigr)
\bigl(M_1^{2j}\bigr)^{{1}/{2}}\bigl(
\Psi^{0,2j}_{n-1}(T) \bigr)^{{1}/{2}}
\\
&&\hspace*{133pt}{}+\bigl\|\sigma'\bigr\|^j_{\infty}
\Psi^{0,j}_{n-1}(T) \biggr).
\end{eqnarray*}
Since $\Psi^{0,2j}_{n-1}(T) \le K(t_{n-1}-s_{n-1})^j$, and $\Psi
^{0,j}_{n-1}(T) \le K(t_{n-1}-s_{n-1})^{{j}/{2}}$, Gronwall's lemma ompletes
the proof.
\end{pf*}

\begin{pf*}{Proof of Remark~\ref{rem6}}
We prove the result for $d=1$. We first prove
that $g(X_T)$ belongs to $ \mathcal{D}^{m,j}$ for all $j\ge2$, that is,
\[
\bigl\|g(X_T)\bigr\|^j_{m,j}=\sum_{l \le m} \sum_{t_1,\ldots,t_l}
\e\bigl[\bigl|D^{(l)}_{t_1,\ldots,t_l} g(X_T)\bigr|^j\bigr]< \infty.
\]
$D^{(l)}_{t_1,\ldots,t_l} g(X_T)$ contains a sum of terms of type
$g^{(k)}(X_T)\prod_{i=1}^k D^{(j_i)}_{t} X_T$, where $k$ varies in
$\{1,\ldots,l\}$,
$|j|_1=l$ and $a(j)=k$ [$a(j)$ denotes the number of nonzero
components of
$j$]. Since $g \in C^{\infty}_p$, and $X$ satisfies \eqref{eq23}, we
get the
result.

Let us now prove that $g(X_T)$ satisfies $\mathcal{H}^2_m$.
$D^{(l_0)}_t\Delta_{t_i,s_i}D^{(l_1)}_s g(X_T)$ contains a sum of terms
of type $g^{(k)}(X_T)\prod_{i=1}^{k-1} (D^{(j_i)}_{t}
X_T) D^{(l'_0)}_t\Delta_{t_i,s_i}D^{(l'_1)}_s X_T$, where $k$ varies in
$\{1,\ldots,l\}$, $|j|_1=l-1-l'_0-l'_1$,
$a(j)=k-1$, $l_0' \le l_0$ and $l_1' \le l_1$. Then, since $g \in
C^{\infty}_p$, $X$ satisfies \eqref{eq23} and \eqref{eq24}, we get $g(X_T)$
satisfies $\mathcal{H}^2_m$, with $\beta_{g(X_T)}=\frac{1}{2}$ and
$k_l^{g(X_T)}$ depends on $(\|g^{(l')}\|_{\infty})_{l' \le l}$, on
$(M^{j'}_{l'})_{l' \le l, j'\le lj}$ and on $K_l^X$.

It remains to prove that $g(X_T)$ satisfies
$\mathcal{H}^3_{p,N}$. $\mathbb{V}(g(X_T))$ is bounded by
$\e((g(X_T))^2)$. Since $g \in
C^{\infty}_p$ and $X$ satisfies $\e(|X_T|^j)<\infty$ for all $j$, we get
that $\mathbb{V}(g(X_T))$ is bounded. We prove that
$\mathbb{V} (g(X_T)\prod_{i=1}^N K_{n_i}(G_i) )$ is bounded
by the
same way.
\end{pf*}

\subsection{Proof of Lemma \protect\texorpdfstring{\ref{lem5}}{4.10}}\label{prooflem5}
We complete the proof for $d=1$.
We prove by induction that $\forall q' \le q$, $(Y^{q'},Z^{q'})$
belongs to $\mathcal{S}^{m,\infty}$,
that is, $\forall j \ge2$
\begin{eqnarray*}
&&\bigl\|\bigl(Y^{q'},Z^{q'}\bigr)\bigr\|^j_{m,j}\\
&&\qquad=
\sum_{1\le l \le m} \sup_{t_1\le\cdots
\le t_l} \biggl\{ \e
\Bigl[\sup_{t_l\le r \le T} \bigl|D^{(l)}_{t_1, \ldots, t_l}
Y^{q'}_r\bigr|^j\Bigr] +\int_{t_l}^T
\e\bigl[\bigl|D^{(l)}_{t_1, \ldots, t_l} Z^{q'}_r\bigr|^j
\bigr] \,dr \biggr\}< \infty.
\end{eqnarray*}
Using \eqref{eqY} gives
\[
D^{(l)}_{t_1, \ldots, t_l} Y^{q'}_r=
\e_r\bigl[ D^{(l)}_{t_1, \ldots, t_l} F^{q'-1}\bigr]-
\int_{t_l}^r D^{(l)}_{t_1, \ldots, t_l} f
\bigl(\theta^{q'-1}_u\bigr) \,du,
\]
where $\theta^{q'-1}_u:=(u,Y^{q'-1}_u,Z^{q'-1}_u
)$.

Using the definition of $F^{q'-1}$ and applying Doob's inequality leads to
\[
\e\Bigl[\sup_{t_l\le r \le T} \bigl|D^{(l)}_{t_1, \ldots, t_l}
Y^{q'}_r\bigr|^j\Bigr] \le C \biggl( \e
\bigl[\bigl|D^{(l)}_{t_1, \ldots, t_l} \xi\bigr|^j\bigr]+\e \biggl(\int
_{t_l}^T \bigl|D^{(l)}_{t_1, \ldots, t_l} f
\bigl(\theta^{q'-1}_u\bigr)\bigr|^j \,du \biggr)
\biggr),
\]
where $C$ is a generic constant depending on $T$ and $j$.

$D^{(l)}_{t_1, \ldots, t_l} f(\theta^{q'-1}_u)$ contains a sum of terms
of type
$\partial^{l_0}_y\,\partial^{l_1}_z f(\theta^{q'-1}_u)\prod_{i=1}^{l_0}
D^{j_i}_t\times\allowbreak  Y_u^{q'-1} \prod_{i=1}^{l_1}
D^{k_i}_t Z_u^{q'-1} $, where $|j|_1+|k|_1=l$, $a(j)=l_0$, $a(k)=l_1$ and
$l_0+l_1 \le l$. Then H\"{o}lder's inequality gives
%
\begin{eqnarray}
\label{eq26} &&\e \biggl(\int_{t_l}^T
\bigl|D^{(l)}_{t_1, \ldots, t_l} f\bigl(\theta^{q'-1}_u
\bigr)\bigr|^j \,du \biggr)
\nonumber
\\[-8pt]
\\[-8pt]
\nonumber
&&\qquad\le C\Biggl(\sum_{k=1}^l
\bigl\|\partial ^k_{\mathrm{sp}}f\bigr\| ^j_{\infty}
\Biggr)\bigl\|\bigl(Y^{q'-1},Z^{q'-1}\bigr)\bigr\|^{lj}_{l,lj}
\end{eqnarray}
and
%
\begin{eqnarray}
\label{eq25} &&\sum_{1\le l \le m} \sup_{t_1\le\cdots\le t_l}
\e\Bigl[\sup_{t_l\le r \le T} \bigl|D^{(l)}_{t_1, \ldots, t_l}
Y^{q'}_r\bigr|^j\Bigr]
\nonumber
\\[-8pt]
\\[-8pt]
\nonumber
&&\qquad\le C \Biggl(\|\xi
\|^j_{m,j}+\sum_{l=1}^m
\Biggl(\sum_{k=1}^l \bigl\|\partial
^k_{\mathrm{sp}}f\bigr\| ^j_{\infty}\Biggr)\bigl \|
\bigl(Y^{q'-1},Z^{q'-1}\bigr)\bigr\|^{lj}_{l,lj}
\Biggr).
\end{eqnarray}

From \eqref{eqY}, we get $D^{(l)}_{t_1, \ldots, t_l}
Z^{q'}_r=\e_r[D^{(l+1)}_{t_1, \ldots, t_l,r} \xi+\int_r^T D^{(l+1)}_{t_1,
\ldots, t_l,r} f(\theta^{q'-1}_u) \,du]$. Then
\begin{eqnarray*}
&&\int_{t_l}^T \e\bigl[\bigl|D^{(l)}_{t_1, \ldots, t_l}
Z^{q'}_r\bigr|^j\bigr] \,dr \\
&&\qquad\le C \biggl(\int
_{t_l}^T \e\bigl[\bigl|D^{(l+1)}_{t_1, \ldots, t_l,r}
\xi\bigr|^j\bigr]\,dr+\int_{t_l}^T \e
\biggl(\biggl\llvert \int_r^T
D^{(l+1)}_{t_1,
\ldots, t_l,r} f\bigl(\theta^{q'-1}_u
\bigr) \,du\biggr\rrvert ^j \biggr)\,dr \biggr).
\end{eqnarray*}
Using \eqref{eq26} yields
\begin{eqnarray*}
&&\sum_{1\le l \le m} \sup_{t_1\le\cdots\le t_l} \int
_{t_l}^T \e \bigl[\bigl|D^{(l)}_{t_1, \ldots, t_l}
Z^{q'}_r\bigr|^j\bigr] \,dr \\
&&\qquad\le C \Biggl( \|\xi
\|^j_{m+1,j}+\sum_{l=1}^m
\Biggl(\sum_{k=1}^{l} \bigl\| \partial
^k_{\mathrm{sp}}f\bigr\|^j_{\infty}\Biggr) \bigl\|
\bigl(Y^{q'-1},Z^{q'-1}\bigr)\bigr\| ^{(l+1)j}_{(l+1),(l+1)j}
\Biggr).
\end{eqnarray*}
Combining this equation with \eqref{eq25} gives
\begin{eqnarray*}
&&\bigl\|\bigl(Y^{q'},Z^{q'}\bigr)\bigr\|^j_{m,j}
\\
&&\qquad\le C\Biggl(\|\xi\|^j_{m+1,j}+\Biggl(\sum
_{k=1}^{m} \bigl\|\partial^k_{\mathrm{sp}}f
\bigr\|^j_{\infty}\Biggr)\sum_{l=1}^m
\bigl\|\bigl(Y^{q'-1},Z^{q'-1}\bigr)\bigr\| ^{(l+1)j}_{(l+1),(l+1)j}
\Biggr).
\end{eqnarray*}
Iterating this inequality yields the result.
We prove that $\forall q' \le q$, $(Y^{q',p},Z^{q',p})$ belongs to
$\mathcal{S}^{m,\infty}$ in the same way. In this case, the generic constant
$C$ depends on $T$, $j$ and $p$, since we need to use the first part of Lemma~\ref{lem4} to upper bound $\e(|\mathcal{C}_{p-l}(D^{(l)}_t F^{(q-1,p)})|^j)$.


\subsection{Proof of Remark \protect\texorpdfstring{\ref{rem5}}{4.13}}\label{proofrem5}
For the sake of clarity, we assume that $\forall i \le r $, $t_{i-1}
\le s_i
\le t_i$ and
$d=1$. Then we show that if $\xi$ satisfies $\mathcal{H}_{r+q}$ and $f
\in C^{0,r+q-1,r+q-1}_b$, then $I_{q,p}:=\int_0^T
f(s,Y^{q,p}_s,Z^{q,p}_s)\,ds$ satisfies
\[
\bigl|\e\bigl(D^{(r)}_{t_1,\ldots,t_r}I_{q,p}\bigr)-\e
\bigl(D^{(r)}_{s_1,\ldots
,s_r}I_{q,p}\bigr)\bigr|\le
K_r^{I_{q,p}}\bigl((t_1-s_1)^{\beta_{I_{q,p}}}+
\cdots+(t_r-s_r)^{\beta
_{I_{q,p}}}\bigr).
\]
Since $I_{0,p}=0$, we deal with the case $q\ge1$.
Since we have
$D^{(r)}_{t_1,\ldots,t_r}I_{q,p}-D^{(r)}_{s_1,\ldots,s_r}I_{q,p}=\sum_{i=1}^r
D^{(i-1)}_{t} \Delta_i D^{(r-i)}_{s} I_{q,p}$, it is enough to prove that
$\e(D^{(i-1)}_{t}\times\allowbreak  \Delta_i D^{(r-i)}_{s} I_{q,p}) \le K_i
(t_i-s_i)^{\beta_{I_{q,p}}}$ (we refer to the beginning of Section~\ref{secappendix} for the definition of $D^{(i-1)}_{t} \Delta_i
D^{(r-i)}_{s}F$).

We introduce $\theta_u^{q,p}=(u,Y_u^{q,p},Z_u^{q,p})$, two vectors $j$
and $m$,
and four integers $k_0$, $k_1$, $l_0$ and $l_1$ such that $l_0\le i-1$,
$l_1 \le r-i$, $|j|_1+|m|_1=r-1-l_0-l_1$ and $k_0+k_1 \le
r$. If $i<r$, $D^{(i-1)}_t \Delta_i D^{(r-i)}_s
I_{q,p}$ contains a sum of terms of type
\[
\int_{s_{r}}^T \partial^{k_0}_y\,
\partial^{k_1}_z f\bigl(\theta^{q,p}_u
\bigr)\prod_{i=1}^{k_0-1} D^{j_i}_{ts}
Y_u^{q,p} \prod_{i=1}^{k_1}
D^{m_i}_{ts} Z_u^{q,p}
\bigl(D^{(l_0)}_t \Delta_i D_s^{(l_1)}
Y^{q,p}_u\bigr) \,du,
\]
where $a(j)=k_0-1$ [$a(j)$ denotes the number of nonzero components of
$j$]
and $a(m)=k_1$ and of type
\[
\int_{s_{r}}^T \partial^{k_0}_y\,
\partial^{k_1}_z f\bigl(\theta^{q,p}_u
\bigr)\prod_{i=1}^{k_0} D^{j_i}_{ts}
Y_u^{q,p} \prod_{i=1}^{k_1-1}
D^{m_i}_{ts} Z_u^{q,p}
\bigl(D^{(l_0)}_t \Delta_i D_s^{(l_1)}
Z^{q,p}_u\bigr) \,du,
\]
where $a(j)=k_0$, $a(m)=k_1-1$. By using the Cauchy--Schwarz
inequality, we get that $\e[D^{(i-1)}_t \Delta_i
D^{(r-i)}_s I_{q,p}] $ is bounded by
\begin{eqnarray*}
\hspace*{-4pt}&&\bigl\|\partial^{k_0+k_1}_{\mathrm{sp}} f\bigr\|_{\infty}\\
\hspace*{-4pt}&&\quad{}\times \e \Biggl(\int
_{s_{r}}^T\prod_{i=1}^{k_0-1}
\bigl(D^{j_i}_{ts} Y_u^{q,p}
\bigr)^2 \prod_{i=1}^{k_1}
\bigl(D^{m_i}_{ts} Z_u^{q,p}
\bigr)^2 \,du \int_{s_{r}}^T
\bigl(D^{(l_0)}_t \Delta_i D_s^{(l_1)}
Y^{q,p}_u\bigr)^2 \,du \Biggr)^{{1}/{2}}
\end{eqnarray*}
(and the same type of term in
$D^{l_0}_t \Delta_i D_s^{(l_1)}
Z^{q,p}_u$) which leads to
%
\begin{eqnarray}
\label{eq29}&&\e\bigl[D^{(i-1)}_t \Delta_i
D^{(r-i)}_s I_{q,p}\bigr]\nonumber\\
&&\qquad \le C\bigl(T, \bigl(\bigl\|
\partial^{k}_{\mathrm{sp}} f\bigr\|_{\infty}\bigr)_{k \le r},
\bigl\|\bigl(Y^{q,p},Z^{q,p}\bigr)\bigr\|_{r-1,2(r-1)}\bigr)
\\
&&\qquad\quad{}\times\sum
_{l_0=0}^{i-1}\sum
_{l_1=0}^{r-i}\sqrt{
\bigl(D^{(l_0)}_t \Delta^{q,p}_i
D_s^{(l_1)} \bigr)_2},\nonumber
\end{eqnarray}
%
where $(D^{(l_0)}_t\Delta_{i}^{q,p} D^{(l_1)}_s)_j:=\e[\sup_{s_r \le u
\le T} |D^{(l_0)}_t \Delta_{i} D_s^{(l_1)} Y^{q,p}_u|^j]
+\e (\int_{s_r}^T |D^{(l_0)}_t\times\allowbreak  \Delta_{i} D_s^{(l_1)}
Z^{q,p}_u|^2\,du )^{{j}/{2}}$.
If $i=r$, $D^{(r-1)}_t \Delta_i I_{q,p}$ contains the same type of integrals
between $s_r$ and $T$ plus an integral between $s_r$ and $t_r$, which is
bounded by $C(T, (\|\partial^{k}_{\mathrm{sp}}
f\|_{\infty})_{k \le r},\|(Y^{q,p},Z^{q,p})\|_{r,2r})(t_r-s_r)$.
Then, since $(Y^{q,p},Z^{q,p}) \in
\mathcal{S}^{r,\infty}$ and $f \in C_b^{0,r+q-1,r+q-1}$, it remains to
take the supremum over
$t_1,\ldots, t_{l_0}, \allowbreak s_{i+1},\ldots,s_{i+l_1}$ in \eqref{eq29} and
to apply
Lemma~\ref{lem8} to end the proof. $K_i$ depends on
$\|\xi\|_{r+q,2{(r+q-1)!}/{(r-1)!}}$,
$(\|\partial_{\mathrm{sp}}^k f\|_{\infty})_{1\le k \le r+q-1}$, $T$ and
$K_r^{q,p}:=\sup_{l\le r} k_l^{q,p}$
(where $k_l^{q,p}$ is defined in Lemma~\ref{lem8}).

\begin{lemme}\label{lem8}
Assume $\xi$ satisfies $\mathcal{H}_{r+q}^2$ and $f \in\mathcal
{C}^{0,r+q-1,r+q-1}_b$.
Then $\forall
i \in\{1,\ldots,r\}$, $\forall l_0 \le i-1$, $\forall l_1 \le r-i$ and
$\forall j \ge2$
\[
\sup_{t_1\le\cdots\le t_{l_0}}\sup_{s_{i+1}\le\cdots\le s_{i+l_1}} \e\bigl[
\bigl(D^{(l_0)}_t \Delta_{i}^{q,p}
D^{(l_1)}_s\bigr)_j\bigr] \le
k_l^{q,p} (t_i-s_i)^{j({1}/{2}\wedge\beta_{\xi})},
\]
where $l=l_0+l_1+1$ and $k^{q,p}_l$ depends on $k^{\xi}_l$,
$T$,$\|\xi\|_{l+q-1,{(l+q-2)!}/{(l-1)!}j}$ and on $(\|\partial
^k_{\mathrm{sp}} f\|_{\infty})_{1\le k \le l+q-2}$.
\end{lemme}

\begin{pf}
We complete the proof by
induction on
$q$. We distinguish\break  cases $l_1>0$ and $l_1=0$. We first consider
$l_1>0$. Let $u$ be in $[s_r,T]$ and\break  $l\le p$ (if~$l>p$, the first term on
the right-hand side of the following equality vanishes). From \eqref{eqYpZp}
and Lemma~\ref{lem4}, we get
$D^{(l_0)}_t \Delta_{i} D^{(l_1)}_s Y^{q,p}_u =\break  \e_u[\mathcal
{C}_{p-r}(D^{(l_0)}_t \Delta_{i} D^{(l_1)}_s F^{q-1,p})]-
\int_{s_{i+l_1}}^u D^{(l_0)}_t \Delta_{i} D^{(l_1)}_s
f(\theta_v^{q-1})\,dv$. Using the definition of $F^{q-1,p}$ [see
\eqref{eqYpZp}], Doob's inequality and Lemma~\ref{lem4} yields
%
\begin{eqnarray}
\label{eq32} &&\e\Bigl[\sup_{u\in[s_r,T]} \bigl(D^{(l_0)}_t
\Delta_{i} D^{(l_1)}_s Y^{q,p}_u
\bigr)^j\Bigr]\nonumber\\
&&\qquad\le C \biggl( \e\bigl[\bigl|D^{(l_0)}_t
\Delta_{i} D^{(l_1)}_s \xi \bigr|^j
\bigr]\\
&&\hspace*{16pt}\qquad\quad{}+\e\biggl[\int_{s_{i+l_1}}^T \bigl|D^{(l_0)}_t
\Delta_{i} D^{(l_1)}_s f\bigl(\theta
^{q-1,p}_v\bigr)\bigr|\,dv\biggr]^j \biggr),\nonumber
\end{eqnarray}
%
where $C$ denotes a generic constant depending on $T$, $j$ and $p$.

Let us now upper bound $\e (\int_{s_r}^T |D^{(l_0)}_t \Delta
_{i} D_s^{(l_1)}
Z^{q,p}_u|^2\,du )^{{j}/{2}}$. Using
\eqref{eqZp} and the
Clark--Ocone formula gives $\int_0^T Z^{q,p}_u
\,dB_u=\cc_p(F^{q-1,p})-\e(\cc_p(F^{q-1,p}))$. Hence, for $v \in
[s_{r},T]$, we have $\int_{s_r}^v Z^{q,p}_u
\,dB_u=\e_{v}(\cc_p(F^{q-1,p}))-\break \e_{s_{r}}(\cc
_p(F^{q-1,p}))=Y^{q,p}_v+\int_{s_{r}}^v
f(\theta^{q-1,p}_u)\,du -Y^{q,p}_{s_{r}}$. Then, we get
\begin{eqnarray*}
\int_{s_r}^v D^{(l_0)}_t
\Delta_{i} D^{(l_1)}_s Z^{q,p}_u
\,dB_u& =& D^{(l_0)}_t \Delta_{i}
D^{(l_1)}_s Y^{q,p}_v-D^{(l_0)}_t
\Delta_{i} D^{(l_1)}_sY^{q,p}_{s_{r}}\\&&{}+
\int_{s_{r}}^v D^{(l_0)}_t
\Delta_{i} D^{(l_1)}_s f\bigl(
\theta_u^{q-1,p}\bigr) \,du.
\end{eqnarray*}
The left-hand side of the Burkh\"{o}lder--Davis--Gundy inequality gives
\begin{eqnarray*}
&&\e \biggl(\int_{s_r}^T \bigl| D^{(l_0)}_t
\Delta_{i} D^{(l_1)}_s Z^{q,p}_u\bigr|^2\,du
\biggr)^{{j}/{2}}\nonumber\\
&&\qquad \le C' \biggl(\e\Bigl[\sup
_{u \in[s_{r},T]
}\bigl |D^{(l_0)}_t \Delta_{i}
D^{(l_1)}_s Y^{q,p}_u\bigl|^j
\Bigr]
\\
&&\hspace*{17pt}\qquad\quad{}+ \e\biggl[\int_{s_r}^T\bigl|D^{(l_0)}_t
\Delta_{i} D^{(l_1)}_s f\bigl(\theta
^{q-1,p}_u\bigr)\bigr|\,du\biggr]^j \biggr),\nonumber
\end{eqnarray*}
where $C'$ denotes a generic constant depending on $T$ and $j$.
Adding $(C'+1)\times$ \eqref{eq32} to the previous equation leads to
%
\begin{eqnarray}
\label{eq33}\quad &&\bigl(D^{(l_0)}_t \Delta_{i}^{q,p}
D^{(l_1)}_s\bigr)_j
\nonumber
\\
&&\qquad\le C \biggl( \e
\bigl[\bigl|D^{(l_0)}_t \Delta_{i} D^{(l_1)}_s
\xi\bigr|^j\bigr]\\
&&\qquad\qquad{}+\e\biggl[\int_{s_{i+l_1}}^T
\bigl|D^{(l_0)}_t \Delta_{i} D^{(l_1)}_s
f\bigl(\theta^{q-1,p}_u\bigr)\bigr|\,du\biggr]^j
\biggr).\nonumber
\end{eqnarray}
We introduce two vectors $j$ and $m$,
and four integers $k_0$, $k_1$, $l'_0$ and $l'_1$\break  such that $l'_0\le l_0$,
$l'_1 \le l_1$, $|j|_1+|m|_1=l-1-l'_0-l'_1$ and $k_0+k_1 \le
l$.\break  $D^{(l_0)}_t \Delta_{i} D^{(l_1)}_s f(\theta^{q-1,p}_u)$ contains a
sum of
terms of type
\[
\partial^{k_0}_y\,\partial^{k_1}_z
f\bigl(\theta^{q-1,p}_u\bigr)\prod
_{i=1}^{k_0-1} D^{j_i}_{ts}
Y_u^{q-1,p} \prod_{i=1}^{k_1}
D^{m_i}_{ts} Z_u^{q-1,p}
\bigl(D^{(l'_0)}_t \Delta_i D_s^{(l'_1)}
Y^{q-1,p}_u\bigr),
\]
where $a(j)=k_0-1$ and $a(m)=k_1$ and of type
\[
\partial^{k_0}_y\,\partial^{k_1}_z
f\bigl(\theta^{q-1,p}_u\bigr)\prod
_{i=1}^{k_0} D^{j_i}_{ts}
Y_u^{q-1,p} \prod_{i=1}^{k_1-1}
D^{m_i}_{ts} Z_u^{q-1,p}
\bigl(D^{(l'_0)}_t \Delta_i D_s^{(l'_1)}
Z^{q-1,p}_u\bigr),
\]
where $a(j)=k_0$, $a(m)=k_1-1$.

By using Cauchy--Schwarz inequality, we get that $\e[\int_{s_{i+l_1}}^T
|D^{(l_0)}_t \Delta_{i} D^{(l_1)}_s \times\allowbreak f(\theta
^{q-1,p}_u)|\,du]^j$ is bounded by
\begin{eqnarray*}
&&\bigl\|\partial^{k_0+k_1}_{\mathrm{sp}} f\bigr\|^j_{\infty}\e
\Biggl(\Biggl(\int_{s_{i+l_1}}^T\prod
_{i=1}^{k_0-1} \bigl(D^{j_i}_{ts}
Y_u^{q-1,p}\bigr)^2 \prod
_{i=1}^{k_1} \bigl(D^{m_i}_{ts}
Z_u^{q-1,p}\bigr)^2 \,du\Biggr)^{{j}/{2}}\\
&&\hspace*{118pt}{}\times \biggl( \int_{s_{l_1+i}}^T \bigl(D^{(l'_0)}_t
\Delta_i D_s^{(l'_1)} Y^{q-1,p}_u
\bigr)^2 \,du\biggr)^{{j}/{2}} \Biggr)
\end{eqnarray*}
(and the same type of term in
$D^{l'_0}_t \Delta_i D_s^{(l'_1)}
Z^{q-1,p}_u$) which leads to
\begin{eqnarray*}
&&\e\biggl[\int_{s_{i+l_1}}^T |D^{(l_0)}_t
\Delta_{i} D^{(l_1)}_s f\bigl(
\theta^{q-1,p}_u\bigr)|\,du\biggr]^j
\\
&&\qquad\le C\bigl(\bigl(\bigl\|\partial^{k}_{\mathrm{sp}} f\bigr\|_{\infty}
\bigr)_{k\le l},\bigl\|\bigl(Y^{q-1,p},Z^{q-1,p}\bigr)
\bigr\|_{l-1,(l-1)j}\bigr) \\
&&\qquad\quad{}\times \sum_{l'_0=0}^{l_0}
\sum_{l'_1=0}^{l_1}\sqrt{
\bigl(D^{(l'_0)}_t \Delta^{q-1,p}_i
D_s^{(l'_1)} \bigr)_{2j}}.
\end{eqnarray*}
It remains to plug this result into \eqref{eq33}, to take the supremum in
$t_1,\ldots,t_{l_0},\break s_{i+1}, \ldots,s_{i+l_1}$ and to apply the induction
hypothesis to obtain
%
\begin{eqnarray}
\label{eq34}\qquad
&&\sup_{t_1\le\cdots\le t_{l_0}}\sup_{s_{i+1}\le\cdots\le s_{i+l_1}} \e
\bigl[\bigl(D^{(l_0)}_t \Delta_{i}^{q,p}
D^{(l_1)}_s\bigr)_j\bigr]
\\
&&\qquad\le
k_{l}^{\xi} (t_i-s_i)^{j\beta_{\xi}}
+C\bigl(\bigl(\bigl\|\partial^{k}_{\mathrm{sp}} f\bigr\|_{\infty}
\bigr)_{1\le k\le l},\bigl\|\bigl(Y^{q-1,p},Z^{q-1,p}\bigr)\bigr\|
_{l-1,(l-1)j}\bigr)
\nonumber
\\[-8pt]
\\[-8pt]
\nonumber
&&\qquad\quad{}\times k_l^{q-1,p}(t_i-s_i)^{j({1}/{2}\wedge
\beta_{\xi})}
\end{eqnarray}
and the result follows. If $l_1=0$, we get
\begin{eqnarray*}
D^{(l_0)}_t \Delta_{i} Y^{q,p}_u
&=& \e_r\bigl[\mathcal{C}_{p-r}\bigl(D^{(l_0)}_t
\Delta_{i} F^{q-1,p}\bigr)\bigr]\\
&&{}-\int_{s_{i}}^u
D^{(l_0)}_t \Delta_{i} D^{(l_1)}_s
f\bigl(\theta_v^{q-1}\bigr)\,dv+\int_{s_i}^{t_i}
D^{(l_0)}_t \Delta_{i} D^{(l_1)}_s
f\bigl(\theta_v^{q-1}\bigr)\,dv.
\end{eqnarray*}
When bounding $\e[\sup_{u \in[s_r, T]}| D^{(l_0)}_t
\Delta_{i} Y^{q,p}_u|^j]$, we deal with the first two terms as we did
before, we bound the term $\e[\int_{s_i}^{t_i}| D^{(l_0)}_t \Delta_i
f(\theta_v^{q-1,p})|\,dv]^j$ by
\[
C\bigl(\bigl(\bigl\|\partial^{k}_{\mathrm{sp}} f\bigr\|_{\infty}
\bigr)_{1 \le k\le l},\bigl\|\bigl(Y^{q-1,p},Z^{q-1,p}\bigr)
\bigr\|_{l,lj}\bigr) (t_i-s_i)^j,
\]
which completes the proof.
\end{pf}

\subsection{Proof of Lemma \protect\texorpdfstring{\ref{lem6}}{4.14}}\label{prooflem6}
We prove the result by induction.
Lemma~\ref{lem6} is true for $p=0$, since $\cc_0^N(F)=\cc_0(F)$.
Assume that
$\e(|(\cc_{p-1}^N -\cc_{p-1})(F)|^2) \le(K^{F}_{p-1})^2
(\frac{T}{N}
)^{2\alpha_F}\sum_{i=1}^{p-1} i^2\frac{T^i}{i!}$. Since we have
\[
\bigl(\cc_p^N -\cc_p\bigr) (F)=\bigl(
\cc_{p-1}^N -\cc_{p-1}\bigr) (F)+
\bigl(P^N_p-P_p\bigr) (F),
\]
it remains to show that $\e(|(P^N_p-P_p)(F)|^2)\le(k^{F}_p)^2
(\frac{T}{N}
)^{2\alpha_F}p^2 \frac{T^p}{p!}$. We recall
%
\begin{equation}\label{eq15}
P_p(F)=\int_0^T\int
_0^{s_p} \cdots\int_0^{s_2}
u_p(s_p,\ldots,s_1)\,dB_{s_1}\cdots
dB_{s_p},
\end{equation}
where $u_p\dvtx s_p,
\ldots,s_1 \longmapsto\e(D^{(p)}_{s_1\cdots
s_p}F
),$
\begin{equation}
\label{eq11}
P^N_p(F)=\sum_{|n|=p}d^n_p
\prod_{1\le i \le N} K_{n_i}(G_i),
\end{equation}
where $d^n_p=n! \e (F \prod
_{1\le i \le N} K_{n_i}(G_i)
).$
Let us rewrite $P^N_p(F)$ as a sum of stochastic integrals.
Let $r \in\nset$. Applying Lemma~\ref{enmart} to $g\dvtx t \longmapsto
\ind_{]\overline{t}_{i-1},\overline{t}_i]}(t)$
yields $M^{r}_t:=h^{r/2}K_{r} (\frac{B_{t}-B_{\overline
{t}_{i-1}}}{\sqrt{h}} )$ is a martingale and
$M^{r}_t=\int_{\overline{t}_{i-1}}^{t} M^{r-1}_s \,dB_s$. Then,
$M^{r}_t=\int_{\overline{t}_{i-1}}^{t} \int_{\overline{t}_{i-1}}^{s_r}
\cdots\int_{\overline{t}_{i-1}}^{s_2} M^{0}_{s_1} \,dB_{s_1}
\cdots dB_{s_r}$. For $r=n_i$ and $t=\overline{t}_i$, we get
\[
K_{n_i}(G_i)=\frac{1}{h^{{n_i}/{2}}}\int_{\overline
{t}_{i-1}}^{\overline{t}_i}
\int_{\overline{t}_{i-1}}^{s_{n_i}} \cdots\int_{\overline
{t}_{i-1}}^{s_2}
\,dB_{s_1} \cdots dB_{s_{n_i}}.
\]
For $|n|:=n_1+\cdots+n_N=p$, we obtain
%
\begin{eqnarray}
\label{eq12} \prod_{1\le i \le N} K_{n_i}(G_i)&=&
\frac{1}{h^{{p}/{2}}}\underbrace{\int_{\overline
{t}_{N-1}}^{T}\cdots
\int_{\overline{t}_{N-1}}^{s_{|n(N-1)|+2}}}_{n_N\
\mathrm{integrals}}\cdots
\nonumber
\\[-8pt]
\\[-8pt]
\nonumber
&&\hspace*{21pt}\underbrace{\int_{\overline{t}_1}^{\overline{t}_2} \cdots\int
_{\overline{t}_1}^{s_{|n(1)|+2}}}_{n_2\ \mathrm{integrals}} \underbrace{\int
_{0}^{\overline{t}_1} \cdots\int_{0}^{s_2}}_{n_1
\ \mathrm{integrals}}
\,dB_{s_1} \cdots dB_{s_{p}},
\\
\label{eq13} d^n_p&=&n!\frac{1}{h^{{p}/{2}}}\underbrace{\int
_{\overline
{t}_{N-1}}^{T}\cdots \int_{\overline{t}_{N-1}}^{l_{|n(N-1)|+2}}}
_{n_N\ \mathrm{integrals}}
\cdots
\nonumber
\\[-8pt]
\\[-8pt]
\nonumber
&&\hspace*{21pt}{} \underbrace{\int_{\overline{t}_1}^{\overline{t}_2} \cdots\int
_{\overline{t}_1}^{l_{|n(1)|+2}}}_{n_2\ \mathrm{integrals}} \underbrace{\int
_{0}^{\overline{t}_1} \cdots\int_{0}^{l_2}}_{n_1\ \mathrm{integrals}}
u_p(l_p,\ldots, l_1) \,dl_1 \cdots
dl_{p}.
\end{eqnarray}
To compare $P_p(F)$ and $P^N_p(F)$, we split the integrals in \eqref{eq15},
%
\begin{eqnarray}
\label{eq14} P_p(F)&=&\sum_{|n|=p}
\underbrace{\int_{\overline{t}_{N-1}}^{T}\cdots \int
_{\overline{t}_{N-1}}^{s_{|n(N-1)|+2}}}_{n_N\
\mathrm{integrals}}\cdots
\nonumber
\\[-8pt]
\\[-8pt]
\nonumber
&&
\hspace*{22pt}\underbrace{\int
_{\overline{t}_1}^{\overline{t}_2} \cdots\int_{\overline{t}_1}^{s_{|n(1)|+2}}}
_{n_2\ \mathrm{integrals}}
\underbrace{\int_{0}^{\overline{t}_1} \cdots\int
_{0}^{s_2}}_{n_1\ \mathrm{integrals}}u_p(s_p,
\ldots,s_1)\,dB_{s_1}\cdots dB_{s_p}.
\end{eqnarray}
Combining \eqref{eq11}, \eqref{eq12}, \eqref{eq13} and \eqref{eq14}
yields $\e(|(P_p^N-P_p)(F)|^2)=$
%
\begin{eqnarray}
\label{eq16} &&\sum_{|n|=p}\underbrace{\int
_{\overline{t}_{N-1}}^{T}\cdots \int_{\overline{t}_{N-1}}^{s_{|n(N-1)|+2}}}_
{n_N \ \mathrm{integrals}}
\cdots
\nonumber
\\[-8pt]
\\[-8pt]
\nonumber
&&\hspace*{22pt}\underbrace{\int_{\overline{t}_1}^{\overline{t}_2} \cdots\int
_{\overline{t}_1}^{s_{|n(1)|+2}}}_{n_2\ \mathrm{integrals}} \underbrace{\int
_{0}^{\overline{t}_1} \cdots\int_{0}^{s_2}}_{n_1
\ \mathrm{integrals}}
\biggl\llvert \frac{d^n_p}{h^{{p}/{2}}}-u_p(s_p,
\ldots,s_1)\biggr\rrvert ^2\,ds_1\cdots
ds_p.
\end{eqnarray}
Moreover, $\frac{d^n_p}{h^{{p}/{2}}}-u_p(s_p,\ldots,s_1)=$
\begin{eqnarray*}
&&\frac{n!}{h^p}\underbrace{\int_{\overline
{t}_{N-1}}^{T}
\cdots \int_{\overline{t}_{N-1}}^{l_{N-1}+1}}_{n_N \ \mathrm{integrals}}\cdots
\\
&&\hspace*{16pt}\underbrace{\int_{\overline{t}_1}^{\overline{t}_2} \cdots\int
_{\overline{t}_1}^{l_{n_1}+1}}_{n_2 \ \mathrm{integrals}} \underbrace{\int
_{0}^{\overline{t}_1} \cdots\int_{0}^{l_2}}_{n_1\ \mathrm{integrals}}
\bigl(u_p(l_p,\ldots,l_1)-u_p(s_p,
\ldots,s_1)\bigr)\,dl_1\cdots dl_p.
\end{eqnarray*}
Since $u_p$ satisfies Hypothesis \ref{hypo3}, we get
$|u_p(l_p,\ldots,l_1)-u_p(s_p,\ldots,s_1)|\le k^F_p (|l_p-s_p|^{\beta_F}
+\cdots+|l_1-s_1|^{\beta_F})\le p k^F_p h^{\beta_F}$. Then
$\llvert \frac{d^n_p}{h^{{p}/{2}}}-u_p(s_p,\ldots,s_1)\rrvert \le
p k^F_p
h^{\beta_F}$. Plugging this result into \eqref{eq16} completes the proof.


\subsection{Proof of Lemma \protect\texorpdfstring{\ref{lem7}}{4.16}}\label{prooflem7}
Using definitions \eqref{chaosdec} and \eqref{chaosdecMC} leads to
\[
\bigl(\cc_p^N-\cc_p^{N,M}\bigr)
(F)=d_0-\hat{d_0}+\sum_{k=1}^p
\sum_{|n|=k} \bigl(d^n_k-
\hat{d^n_k}\bigr)\prod_{i=1}^N
K_{n_i}(G_i).
\]
Since $\hat{d^n_k}$ is independent of $(G_i)_i$,
\[
\e\bigl(\bigl|\bigl(\cc_p^N-\cc_p^{N,M}
\bigr) (F)\bigr|^2\bigr)=\e\bigl(|d_0-\hat{d_0}|^2
\bigr)+\sum_{k=1}^p \sum
_{|n|=k}\frac{1}{n!}\e\bigl(\bigl|d^n_k-
\hat{d^n_k}\bigr|^2\bigr).
\]
The definition of the coefficients $d_0$ and $d^n_k$ given in
\eqref{coefchaosdec} leads to
\[
\e\bigl(\bigl|\bigl(\cc_p^N-\cc_p^{N,M}
\bigr) (F)\bigr|^2\bigr)=\V(\hat{d_0})+\sum
_{k=1}^p \sum_{|n|=k}
\frac{1}{n!}\V\bigl(\hat{d^n_k}\bigr),
\]
and the first result follows.
To get the second result, we write
$\cc_p^{N,M}(F)=(\cc_p^{N,M}-\cc_p^{N})(F)+\cc_p^N(F)$. Since
$\e ((\cc_p^{N,M}-\cc_p^{N})(F)\cc_p^N(F) )=0$, we get
\[
\e\bigl(\bigl|\cc_p^{N,M}(F)\bigr|^2\bigr)=\e\bigl(\bigl|
\bigl(\cc_p^{N,M}-\cc_p^{N}\bigr)
(F)\bigr|^2\bigr)+\e\bigl(\bigl|\cc _p^N(F)\bigr|^2
\bigr).
\]
Lemma~\ref{lem9} completes the proof.

\section{Wiener chaos expansion formulas}
\subsection{Proof of Proposition \protect\texorpdfstring{\ref{prop2}}{2.7}} 
\label{subbmR}

First, we compute $\e_t(\cc^N_p(F))$ for $t\in\allowbreak ]\overline{t}_{r-1},
\overline{t}_r]$.
From \eqref{chaosdec}, we get
\[
\e_t\bigl(\cc_p^N F\bigr)  =
d_0 + \sum_{k=1}^p \sum
_{|n|=k} d_k^n \prod
_{i<r} K_{n_i}(G_i)\times
\e_t \biggl(\prod_{i\geq r}
K_{n_i}(G_i) \biggr).
\]
Since Brownian
increments are independent, we get $\e_{\overline{t}_r}(\prod_{i\geq r}
K_{n_i}(G_i)) =\allowbreak  K_{n_r}(G_r) \prod_{i>r}
\e[K_{n_i}(G_i)]$, which is null as soon as $n_{r+1} + \cdots+ n_N
>0$. Then, nested conditional expectations give
\[
\e_t\bigl(\cc_p^N F\bigr)  =
d_0 + \sum_{k=1}^p \sum
_{|n(r)|=k} d_k^n \prod
_{i<r} K_{n_i}(G_i)\times
\e_t \bigl( K_{n_r}(G_r) \bigr).
\]
By applying Lemma~\ref{enmart} when $g\dvtx t \longmapsto\ind
_{]\overline
{t}_{r-1},\overline{t}_r]}(t)$, we get $\e_t  (
K_{n_r}(G_r) ) =  (\frac{t-\overline{t}_{r-1}}{h} )^{n_r/2}
K_{n_r} (\frac{B_t-B_{\overline{t}_{r-1}}}{\sqrt{t-\overline
{t}_{r-1}}} )$, which yields
the first result.
Since $K'_n(x) =\break  K_{n-1}(x)$, the second result follows.


\subsection{Wiener chaos expansion formulas in $\rset^d$} 
\label{subbmRd}

We want to approximate $F\in\lp^2(\m F_T)$ using its chaos
decomposition up to
order $p$. We assume $N\geq dp$. We consider the following truncated
basis of $\lp^2 ([0,T]; \rset^d )$:
\[
\frac{\ind_{]\overline{t}_{i-1}, \overline{t}_i]}(t)}{\sqrt{h}} e_j,\qquad i=1,\ldots,N, j=1,\ldots,d, \mbox{ where }
h=\frac{T}{N},
\]
where $\{\overline{t}_i:=ih, i=0,\ldots,N\}$ is a regular mesh grid,
and $(e_j)_{1\leq j\leq d}$ represents the canonical basis of
$\rset^d$. $P_k$, the $k$th chaos, is generated by
\[
\Biggl\{ \prod_{j=1}^d \prod
_{i=1}^N K_{n_i^j} \bigl(G_i^j
\bigr) \dvtx \sum_{j=1}^d\sum
_{i=1}^N n_i^j = k
\Biggr\},\qquad G_i^j = \frac{\Delta
_i^j}{\sqrt{h}},
\Delta_i^j =B^j_{\overline{t}_i}-B^j_{\overline{t}_{i-1}}.
\]
For $j=1,\ldots,d$, $n^j=(n^j_1,\ldots,n^j_N)$, one notes
$|n^j|=n^j_1+\cdots+n^j_N$, $n^j!=n^j_1!\cdots n^j_N!$, and for $r\leq N$,
$n^j(r)=(n^j_1,\ldots,n^j_r)$. $n=(n^1,\ldots,n^d)^*$, $|n| =
|n^1|+\cdots+|n^d|$, $n!=n^1!\cdots n^d!$ and
$n(r)=(n^1(r),\ldots,n^d(r))^*$. Since the r.v.\break  $(\prod_{1\leq
j\leq d}
\prod_{1\leq i\leq N} K_{n_i^j} (G_i^j ) )_n$ are orthogonal
ones, the projection of $F$ is given by
\[
\cc^N_p (F) = d_0 + \sum
_{k=1}^p \sum_{|n|=k}
d_k^n \prod_{1\leq
j\leq d}\prod
_{1\leq i\leq N} K_{n_i^j} \bigl(G_i^j
\bigr),
\]
where the coefficients $d_k^n$ are given by
\[
d_k^n = n! \e \biggl[F \prod
_{1\leq j\leq d}\prod_{1\leq
i\leq N}
K_{n_i^j} \bigl(G_i^j \bigr) \biggr].
\]

\begin{prop}\label{prop3}
For $\overline{t}_{r-1}<t\leq\overline{t}_r$, we have
\begin{eqnarray*}
\e_t\bigl(\cc^N_p F\bigr) &= &d_0
+ \sum_{k=1}^p \sum
_{|n(r)|=k} d_k^n \prod
_{i<r} \prod_{1\leq j\leq d}K_{n_i^j}
\bigl(G_i^j \bigr)\\
&&\hspace*{19pt}{}\times\prod
_{1\leq
j\leq d} \biggl(\frac{t-\overline{t}_{r-1}}{h} \biggr)^{{n_r^j}/{2}}
K_{n_r^j} \biggl(\frac{B^j_t-B^j_{\overline{t}_{r-1}}}{\sqrt {t-\overline
{t}_{r-1}}} \biggr),
\end{eqnarray*}
and for $l=1,\ldots,d$,
\begin{eqnarray*}
&&D_t^l\bigl(\e_t\bigl(\cc^N_p
F\bigr)\bigr)\\
&&\qquad = \sum_{k=1}^p \mathop{
\sum_{|n(r)|=k }}_{
n_r^l>0} d_k^n
h^{-1/2} \prod_{i<r}\prod
_{1\leq j\leq d} K_{n_i^j} \bigl(G_i^j
\bigr)
 \biggl(\frac{t-\overline{t}_{r-1}}{h} \biggr)^{
{(n_r^l-1)}/{2}} \\
&&\qquad\quad{}\times K_{n_r^l-1}
\biggl(\frac{B^l_t-B^l_{\overline
{t}_{r-1}}}{\sqrt{t-\overline{t}_{r-1}}} \biggr)\\
&&\qquad\quad{}\times \prod_{j\neq l} \biggl(
\frac{t-\overline{t}_{r-1}}{h} \biggr)^{{n_r^j}/{2}} K_{n_r^j} \biggl(
\frac{B^j_t-B^j_{\overline{t}_{r-1}}}{\sqrt{t-\overline
{t}_{r-1}}} \biggr).
\end{eqnarray*}
\end{prop}

\begin{rem}
In particular, for $t=\overline{t}_r$, $r \ge1$ and $l=1,\ldots,d$,
\begin{eqnarray*}
\e_{\overline{t}_r}\bigl(\cc^N_p F\bigr) & =&
d_0 + \sum_{k=1}^p \sum
_{|n(r)|=k} d_k^n \prod
_{i\leq r} \prod_{1\leq j\leq d}K_{n_i^j}
\bigl(G_i^j \bigr),
\\
D^l_{\overline{t}_r}\bigl(\e_{\overline{t}_r} \bigl(
\cc^N_p F\bigr)\bigr) & = &\sum
_{k=1}^p \mathop{\sum_{|n(r)|=k }}_{ n_r^l>0}
d_k^n h^{-1/2} \prod
_{i<r}\prod_{1\leq j\leq d}
K_{n_i^j} \bigl(G_i^j \bigr)
K_{n_r^l-1} \bigl(G^l_r \bigr) \prod
_{j\neq l} K_{n_r^j} \bigl(G^j_r
\bigr).
\end{eqnarray*}
When $r=0$, we get $\e_{\overline{t}_0}(\cc^N_p F)=d_0$, and we
define $
D^l_{\overline{t}_0}(\e_{\overline{t}_0} (\cc^N_p
F))=\frac{1}{\sqrt{h}}d_1^{e^l_1}$, where $(e^i_j)$ is a matrix of size
$d\times N$ whose component $(i,j)$ equals $1$ and the other ones are null.
\end{rem}

\begin{pf*}{Proof of Proposition~\ref{prop3}}
We first compute $\e_t(\cc^N_p F)$ for $t\in\,]\overline{t}_{r-1},
\overline{t}_r]$. We have
\[
\e_t\bigl(\cc^N_p F\bigr) =
d_0 + \sum_{k=1}^p \sum
_{|n|=k} d_k^n \prod
_{i<r} \prod_{1\leq j\leq d}K_{n_i^j}
\bigl(G_i^j \bigr)\times\e_t \biggl(\prod
_{i\geq r} \prod_{1\leq j\leq d}G_{n_i^j}
\bigl(W_i^j \bigr) \biggr).
\]
Since Brownian motions and their increments are independents, we get
\[
\e_{\overline{t}_r} \biggl(\prod_{i\geq r} \prod
_{1\leq
j\leq
d}K_{n_i^j} \bigl(G_i^j
\bigr) \biggr) = \prod_{1\leq j\leq
d}K_{n_r^j}
\bigl(G_r^j \bigr) \prod_{i>r}
\prod_{1\leq j
\leq d} \e \bigl[K_{n_i^j}
\bigl(G_i^j \bigr) \bigr],
\]
which is null as soon as $n_{r+1}^1 + \cdots+ n_N^1 + \cdots+ n_{r+1}^d +
\cdots+ n_N^d>0$. Then nested conditional expectations give
\[
\e_t(F)  = d_0 + \sum_{k=1}^p
\sum_{|n(r)|=k} d_k^n \prod
_{i<r} \prod_{1\leq j\leq d}
K_{n_i^j} \bigl(G_i^j \bigr)\times
\e_t \biggl( \prod_{1\leq j\leq d}K_{n_r^j}
\bigl(G_r^j \bigr) \biggr).
\]
From Lemma~\ref{enmart}, for $j=1,\ldots,d$ $M^{n_r^j}_t:=
(t-\overline{t}_{r-1} )^{n_r^j/2}
K_{n_r^j} (\frac{B^j_t-B^j_{\overline{t}_{r-1}}}{\sqrt {t-\overline
{t}_{r-1}}} )$ is a
martingale and $ dM^{n_r^j}_t = M^{n_r^j-1}_t \ind_{]\overline
{t}_{r-1},\overline{t}_r]}(t)
\,dB^j_t$. Then $\prod_{1\leq j\leq d}  (t-\overline
{t}_{r-1}
)^{n_r^j/2}\times\allowbreak
K_{n_r^j} (\frac{B^j_t-B^j_{\overline{t}_{r-1}}}{\sqrt {t-\overline
{t}_{r-1}}} )$\vspace*{2pt} is also a
martingale, and the first result follows. Since\break  $K'_{n_r^l}(x) =
K_{n_r^l-1}(x)$, we get the second result.
\end{pf*}
\end{appendix}


%

%



\printaddresses

\end{document}